\documentclass[11pt,a4paper,reqno]{amsart}%
\usepackage{amsfonts}
\usepackage{amsmath}
\usepackage{amssymb}
\usepackage{graphicx}%
\newtheorem{theorem}{Theorem}
\theoremstyle{plain}

\newtheorem{corollary}{Corollary}

\newtheorem{definition}{Definition}
\newtheorem{example}{Example}

\newtheorem{lemma}{Lemma}

\newtheorem{remark}{Remark}

\numberwithin{equation}{section}
\begin{document}
\title{Shift operators and stability in delayed dynamic equations }
\author{Murat Ad\i var}
\address{Izmir University of Economics\\
Department of Mathematics, 35330, Izmir Turkey}
\email[M. Ad\i var]{murat.adivar@ieu.edu.tr}
\author{Youssef N. Raffoul}
\address{University of Dayton\\
Department of Mathematics, Dayton, OH 45469-2316, USA}
\email[Y. Raffoul]{ youssef.raffoul@notes.udayton.edu}

\begin{abstract}
In this paper, we use what we call the shift operator so that general delay
dynamic equations of the form
\[
x^{\Delta}(t)=a(t)x(t)+b(t)x(\delta_{-}(h,t))\delta_{-}^{\Delta}%
(h,t),\ \ \ t\in\lbrack t_{0},\infty)_{\mathbb{T}}%
\]
can be analyzed with respect to stability and existence of solutions. By means
of the shift operators we define a general delay function opening an avenue
for the construction of Lyapunov functional on time scales. Thus, we use the
Lyapunov's direct method to obtain inequalities that lead to stability and
instability. Therefore, we extend and unify stability analysis of delay
differential, delay difference, delay $h-$difference, and delay $q-$difference
equations which are the most important particular cases of our delay dynamic equation.

\end{abstract}
\maketitle

\section{Introduction}

Lyapunov functionals are widely used in stability analysis of differential and
difference equations. However, the extension of utilization of Lyapunov
functionals in dynamical systems on time scales has been lacking behind due to
the constrained presented by the particular time scale. For example, in delay
differential equations, a suitable Lyapunov functional will involve a term
with double integrals, in which one of the integral's lower limit is of the
form $t+s$. Such a requirement will restrict the time scale that can be considered.

For a few references on the study of stability in differential equations,
using Lyapunov functionals, we refer the interested reader to \cite{mayr},
\cite{raffoul&adivar}, \cite{burton 1}-\cite{yoshizawa}. The reader may
consult Yoshizawa \cite[pp. 183-213]{yoshizawa} (or any book on functional
differential equations and Lyapunov's direct method) for definitions of
stability and for properties of Lyapunov functionals. For the stability
analysis of the delay differential equation%
\begin{equation}
x^{\prime}(t)=a(t)x(t)+b(t)x(t-h),\ \ h>0\label{delay differential}%
\end{equation}
we refer to \cite{busenberg}-\cite{hatvani}, and \cite{wang}. In
\cite{adraf2}, the authors improved the results of \cite{wang} by considering
the delay differential equation of the form%
\begin{equation}
x^{\prime}(t)=a(t)x(t)+b(t)x(t-h(t)),\ \ 0<h(t)\leq r_{0}%
.\label{delay diff time}%
\end{equation}
On the other hand, stability analysis of delay difference equations of the
form
\begin{equation}
x(t+1)=a(t)x(t)+b(t)x(t-\tau),\ \ \tau\in\mathbb{Z}_{+}%
\label{delay difference}%
\end{equation}
is treated in \cite{berezansky}, \cite{raffoul}, and \cite{raffoul2}.

A time scale, denoted $\mathbb{T}$, is a nonempty closed subset of real
numbers. The set $\mathbb{T}^{\kappa}$ is derived from the time scale
$\mathbb{T}$ as follows: if $\mathbb{T}$ has a left-scattered maximum $M$,
then $\mathbb{T}^{\kappa}=\mathbb{T-}\left\{  M\right\}  $, otherwise
$\mathbb{T}^{\kappa}=\mathbb{T}$. The delta derivative $f^{\Delta}$ of a
function $f:\mathbb{T\rightarrow R}$, defined at a point $t\in\mathbb{T}%
^{\kappa}$ by%
\begin{equation}
f^{\Delta}(t):=\lim_{s\rightarrow t}\frac{f(\sigma(t))-f(s)}{\sigma
(t)-s}\text{,\ \ \ where }s\rightarrow t\text{,\ \ }s\in\mathbb{T}%
\backslash\left\{  \sigma(t)\right\}  , \label{delta derivative}%
\end{equation}
was first introduced by Hilger \cite{hilger} to unify discrete and continuous
analyses. In (\ref{delta derivative}), $\sigma:\mathbb{T}\rightarrow
\mathbb{T}$ is the forward jump operator defined by $\sigma(t):=\inf\left\{
s\in\mathbb{T}:s>t\right\}  $. Hereafter, we denote by $\mu(t)$ the step size
function $\mu:\mathbb{T}\rightarrow\mathbb{R}$ defined by $\mu(t):=\sigma
(t)-t$. A point $t\in\mathbb{T}$ is said to be right dense (right scattered)
if $\mu(t)=0$ ($\mu(t)>0$). A point is said to be left dense if $\sup\left\{
s\in\mathbb{T}:s<t\right\}  =t$. A function $f:\mathbb{T}\rightarrow
\mathbb{R}$ is called $rd$-continuous if it is continuous at right dense
points and its left sided limits exists (finite) at left dense points. Every
$rd$-continuous function $f:\mathbb{T}\rightarrow\mathbb{R}$ has an
anti-derivative $F$ denoted by%
\[
F(t):=\int_{t_{0}}^{t}f(t)\Delta t.
\]
To indicate the time scale interval $[a,b]\cap\mathbb{T}$ we use the notation
$\left[  a,b\right]  _{\mathbb{T}}$. The intervals $[a,b)_{\mathbb{T}}$,
$(a,b]_{\mathbb{T}}$, and $\left(  a,b\right)  _{\mathbb{T}}$ are defined
similarly. For brevity, we assume the reader is familiar with the basic
calculus of time scales. \ A comprehensive review on dynamic equations on time
scales can be found in \cite{book} and \cite{book2}.

In \cite{mayr} and \cite{Al}, the authors handle the stability analysis of the
dynamic equation%
\begin{equation}
x^{\Delta}(t)=a(t)x(t)+b(t)x(\delta(t))\delta^{\Delta}(t),
\label{delay dynamic}%
\end{equation}
where the delay function $\delta:[t_{0},\infty)_{\mathbb{T}}\rightarrow
\lbrack\delta(t_{0}),\infty)_{\mathbb{T}}$ is surjective, strictly increasing
and is supposed to have the following properties%
\[
\delta(t)<t,\ \ \ \ \delta^{\Delta}(t)<\infty,\ \ \ \delta\circ\sigma
=\sigma\circ\delta\text{.}%
\]
Afterwards, we point out in \cite{mayr2} that the assumption $\delta
\circ\sigma=\sigma\circ\delta$ is redundant whenever the delay function
$\delta:[t_{0},\infty)_{\mathbb{T}}\rightarrow\lbrack\delta(t_{0}%
),\infty)_{\mathbb{T}}$ is surjective and strictly increasing.

Note that the delta derivative in (\ref{delta derivative}) turns into the
ordinary derivative $f^{\prime}(t)$ and the forward difference $\Delta
f(t):=f(t+1)-f(t)$ when $\mathbb{T}=\mathbb{R}$ and $\mathbb{T}=\mathbb{Z}$,
respectively. Hence, (\ref{delay dynamic}) is a general equation including the
particular cases (\ref{delay differential})-(\ref{delay difference}). However,
this paper improves the results of \cite{mayr}.

In this paper, we define the general shift operator and make use of them in
the construction of the Lyapunov functional to improve previous results on
delay dynamic equations regarding stability and boundedness of solutions. In
particular, we improve the results of Eq. (\ref{delay differential}%
)-(\ref{delay difference}), and \ref{delay dynamic}). The main task of this
paper can be outlined as follows:

\begin{itemize}
\item To create a suitable Lyapunov function that leads to exponential
stability of the zero solution.

\item To give criteria for instability.

\item To compare the results of this paper with the ones in the existing literature.
\end{itemize}

in \cite{wang}, the author used the following
\begin{align}
V(t)  &  =\left[  x(t)+\int_{t-h}^{t}b(s+h)x(s)ds\right]  ^{2}\nonumber\\
&  +\lambda\int_{-h}^{0}\int_{t+s}^{t}b^{2}(z+h)x^{2}(z)dzds.
\label{lyapunov continuous1}%
\end{align}
to study the exponential stability of the zero solution of \eqref{delay differential}.

We do not adopt this type of Lyapunov functional since it requires the time
scale to be additive. An additive time scale is a time scale which is closed
under addition. There are many time scales that are not additive. To be more
specific, the time scales $\overline{q^{Z}}=\left\{  0\right\}  \cup\left\{
q^{n}:n\in\mathbb{Z}\right\}  $,\ $\sqrt{\mathbb{N}}=\left\{  \sqrt{n}%
:n\in\mathbb{N}\right\}  $ are not additive. However, $\delta_{\pm
}(s,t)=ts^{\pm1}$ and $\delta_{\pm}(s,t)=\sqrt{t^{2}\pm s^{2}}$ are the shift
operators defined on $\overline{q^{Z}}$ and $\sqrt{\mathbb{N}}$, respectively.
It turns out that we need the notion of shift operators to avoid additivity
assumption on the time scale. That is, to include more time scales in the
investigation. Shift operators are first introduced in \cite{adivar} to obtain
function bounds for convolution type Volterra integro-dynamic equations on
time scales. However, the time scales considered in \cite{adivar} is
restricted to the ones having an initial point $t_{0}\in\mathbb{T}$ so that
there exist the shift operators defined on $[t_{0},\infty)\cap\mathbb{T}$.
Afterwards, in \cite{bams} the definition of shift operators was extended so
that they are defined on the whole time scale $\mathbb{T}$. In this paper, our
new and generalized shift operators include positive and negative values.

We end this section by giving some basic definitions and theorems that will be
used in further sections.

\begin{definition}
A function $h:\mathbb{T}\rightarrow\mathbb{R}$ is said to be \emph{regressive}
provided $1+\mu(t)h(t)\neq0$ for all $t\in\mathbb{T}^{\kappa}$, where
$\mu(t)=\sigma(t)-t$. The set of all regressive $rd$-continuous functions
$\varphi:\mathbb{T}\rightarrow\mathbb{R}$ is denoted by $\mathcal{R}$ while
the set $\mathcal{R}^{+}$ is given by $\mathcal{R}^{+}=\{h\in\mathcal{R}%
:1+\mu(t)\varphi(t)>0\mbox{
for all }t\in\mathbb{T}\}$.
\end{definition}

Let $\varphi\in\mathcal{R}$ and $\mu(t)>0$ for all $t\in\mathbb{T}$. The
\emph{exponential function} on $\mathbb{T}$ is defined by
\begin{equation}
e_{\varphi}(t,s)=\exp\left(  \int_{s}^{t}\!\zeta_{\mu(r)}(\varphi(r))\Delta
r\right)  \label{exp}%
\end{equation}
where $\zeta_{\mu(s)}$ is the cylinder transformation given by
\begin{equation}
\zeta_{\mu(r)}(\varphi(r))\!:=\left\{
\begin{array}
[c]{cc}%
\frac{1}{\mu(r)}\mbox{Log}(1+\mu(r)\varphi(r)) & if\text{ }\mu(r)>0\\
\varphi(r) & if\ \ \mu(r)=0
\end{array}
\right.  \,.\label{cylinder}%
\end{equation}
It is well known that if $p\in\mathcal{R}^{+}$, then $e_{p}(t,s)>0$ for all
$t\in\mathbb{T}$. Also, the exponential function $y(t)=e_{p}(t,s)$ is the
solution to the initial value problem $y^{\Delta}=p(t)y,\,y(s)=1$. Other
properties of the exponential function are given in the following lemma:

\begin{lemma}
\label{lemma2.3} \cite[Theorem 2.36]{book} Let $p,q\in\mathcal{R}$. Then

\begin{itemize}
\item[i.] $e_{0}(t,s)\equiv1$ and $e_{p}(t,t)\equiv1$;

\item[ii.] $e_{p}(\sigma(t),s)=(1+\mu(t)p(t))e_{p}(t,s)$;

\item[iii.] $\frac{1}{e_{p}(t,s)}=e_{\ominus p}(t,s)$ where, $\ominus
p(t)=-\frac{p(t)}{1+\mu(t)p(t)}$;

\item[iv.] $e_{p}(t,s)=\frac{1}{e_{p}(s,t)}=e_{\ominus p}(s,t)$;

\item[v.] $e_{p}(t,s)e_{p}(s,r)=e_{p}(t,r)$;

\item[vi.] $\left(  \frac{1}{e_{p}(\cdot,s)}\right)  ^{\Delta}=-\frac
{p(t)}{e_{p}^{\sigma}(\cdot,s)}$.
\end{itemize}
\end{lemma}

\begin{theorem}
\cite[Theorem 1.117]{book} \label{theorem 1.117} Let $a\in\mathbb{T}^{\kappa}%
$, $b\in\mathbb{T}$ and assume that $k:\mathbb{T}\times\mathbb{T}^{\kappa
}\rightarrow\mathbb{R}$ is continuous at $\left(  t,t\right)  $, where
$t\in\mathbb{T}^{\kappa}$ with $t>a$. Also assume that $k^{\Delta}\left(
t,.\right)  $ is rd-continuous on $\left[  a,\sigma\left(  t\right)  \right]
$. Suppose that for each $\varepsilon>0$ there exists a neighborhood $U$ of
$t$, independent of $\tau\in\left[  t_{0},\sigma\left(  t\right)  \right]  $,
such that
\[
\left\vert k\left(  \sigma\left(  t\right)  ,\tau\right)  -k\left(
s,r\right)  -k^{\Delta}\left(  t,\tau\right)  \left(  \sigma\left(  t\right)
-s\right)  \right\vert \leq\varepsilon\left\vert \sigma\left(  t\right)
-s\right\vert
\]
for all $s\in U$, where $k^{\Delta}$ denotes the derivative of $k$ with
respect to the first variable. Then
\begin{align*}
g\left(  t\right)   &  :=\int_{a}^{t}k\left(  t,\tau\right)  \Delta\tau\text{
implies }g^{\Delta}\left(  t\right)  =\int_{a}^{t}k^{\Delta}\left(
t,\tau\right)  \Delta\tau+k\left(  \sigma\left(  t\right)  ,t\right) \\
h\left(  t\right)   &  :=\int_{t}^{b}k\left(  t,\tau\right)  \Delta\tau\text{
implies }g^{\Delta}\left(  t\right)  =\int_{t}^{b}k^{\Delta}\left(
t,\tau\right)  \Delta\tau-k\left(  \sigma\left(  t\right)  ,t\right)  .
\end{align*}

\end{theorem}

\section{Shift operators}

Next, we state the generalized shift operators. A limited version of it can be
found in \cite{adivar}.

\begin{definition}
\label{shift} Let $\mathbb{T}^{\ast}$ be a non-empty subset of the time scale
$\mathbb{T}$ and $t_{0}\in\mathbb{T}^{\ast}$ a fixed number such that there
exist operators $\delta_{\pm}:[t_{0},\infty)_{\mathbb{T}}\times\mathbb{T}%
^{\ast}\rightarrow\mathbb{T}^{\ast}$ satisfying the following properties:

\begin{enumerate}
\item[P.1] The functions $\delta_{\pm}$ are strictly increasing with respect
to their second arguments, i.e., if
\[
(T_{0},t),(T_{0},u)\in\mathcal{D}_{\pm}:=\left\{  (s,t)\in\lbrack t_{0}%
,\infty)_{\mathbb{T}}\times\mathbb{T}^{\ast}:\delta_{\pm}(s,t)\in
\mathbb{T}^{\ast}\right\}  ,
\]
then
\[
T_{0}\leq t<u\text{ implies }\delta_{\pm}(T_{0},t)<\delta_{\pm}(T_{0},u),
\]

\item[P.2] If $(T_{1},u),(T_{2},u)\in\mathcal{D}_{-}$ with $T_{1}<T_{2}$, then%
\[
\delta_{-}(T_{1},u)>\delta_{-}(T_{2},u),
\]
and if $(T_{1},u),(T_{2},u)\in\mathcal{D}_{+}$ with $T_{1}<T_{2}$, then
\[
\delta_{+}(T_{1},u)<\delta_{+}(T_{2},u),
\]

\item[P.3] If $t\in\lbrack t_{0},\infty)_{\mathbb{T}}$, then $(t,t_{0}%
)\in\mathcal{D}_{+}$ and $\delta_{+}(t,t_{0})=t$. Moreover, if $t\in
\mathbb{T}^{\ast}$, then $(t_{0},t)$ $\in\mathcal{D}_{+}$ and $\delta
_{+}(t_{0},t)=t$ holds,

\item[P.4] If $(s,t)\in\mathcal{D}_{\pm}$, then $(s,\delta_{\pm}%
(s,t))\in\mathcal{D}_{\mp}$ and $\delta_{\mp}(s,\delta_{\pm}(s,t))=t$,

\item[P.5] If $(s,t)\in\mathcal{D}_{\pm}$ and $(u,\delta_{\pm}(s,t))\in
\mathcal{D}_{\mp}$, then $(s,\delta_{\mp}(u,t))\in\mathcal{D}_{\pm}$ and
\[
\delta_{\mp}(u,\delta_{\pm}(s,t))=\delta_{\pm}(s,\delta_{\mp}(u,t)).
\]

\end{enumerate}

\noindent Then the operators $\delta_{-}$ and $\delta_{+}$ associated with
$t_{0}\in\mathbb{T}^{\ast}$ (called the initial point) are said to be
\textit{backward and forward shift operators} on the set $\mathbb{T}^{\ast}$,
respectively. The variable $s\in\lbrack t_{0},\infty)_{\mathbb{T}}$ in
$\delta_{\pm}(s,t)$ is called the shift size. The values $\delta_{+}(s,t)$ and
$\delta_{-}(s,t)$ in $\mathbb{T}^{\ast}$ indicate $s$ units translation of the
term $t\in\mathbb{T}^{\ast}$ to the right and left, respectively. The sets
$\mathcal{D}_{\pm}$ are the domains of the shift operators $\delta_{\pm}$, respectively.
\end{definition}

\begin{definition}
Let $\mathbb{T}$ be a time scale having an initial point such that there exist
operators $\delta_{\pm}:[t_{0},\infty)_{\mathbb{T}}\times\mathbb{T}%
\rightarrow\mathbb{T}$ satisfying P.3-P.5. A point $t^{\ast}(\neq t_{0}%
)\in\mathbb{T}$ is said to be a sticky point of $\mathbb{T}$ if%
\[
\delta_{\pm}(s,t^{\ast})=t^{\ast}\text{ for all }s\in\lbrack t_{0}%
,\infty)_{\mathbb{T}}\text{ with }(s,t^{\ast})\in\mathcal{D}_{\pm}.
\]

\end{definition}

Hereafter, let $t^{\ast}$ and $\mathbb{T}^{\ast}$ denote the sticky point and
the largest subset of $\mathbb{T}$ without sticky point, respectively.

\begin{corollary}
A sticky point $t^{\ast}$ cannot be included in the interval $[t_{0}%
,\infty)_{\mathbb{T}}$.
\end{corollary}

\begin{proof}
First, by P.3-P.5 we have $\delta_{-}(u,u)=\delta_{-}(u,\delta_{+}%
(u,t_{0}))=t_{0}$, and hence, $(u,u)\in\mathcal{D}_{-}$ for all $u\in\lbrack
t_{0},\infty)_{\mathbb{T}}$. If $t^{\ast}\in\lbrack t_{0},\infty)_{\mathbb{T}%
}$ is a sticky point, then P.3-P.5 imply%
\[
t^{\ast}=\delta_{-}(t^{\ast},t^{\ast})=t_{0}\in\mathbb{T}^{\ast}%
=\mathbb{T-}\left\{  t^{\ast}\right\}  .
\]
This leads to a contradiction.
\end{proof}

\begin{example}
Let $\mathbb{T=R}$ and $t_{0}=1$. The operators%
\begin{equation}
\delta_{-}(s,t)=\left\{
\begin{array}
[c]{cc}%
t/s & \text{if }t\geq0\\
st & \text{if }t<0
\end{array}
\right.  ,\ \ \ \text{for }s\in\lbrack1,\infty) \label{rs1}%
\end{equation}
and%
\begin{equation}
\delta_{+}(s,t)=\left\{
\begin{array}
[c]{cc}%
st & \text{if }t\geq0\\
t/s & \text{if }t<0
\end{array}
\right.  ,\ \ \ \text{for }s\in\lbrack1,\infty) \label{rs2}%
\end{equation}
are backward and forward shift operators associated with the initial point
$t_{0}=1$. Also, $t^{\ast}=0$ is a sticky point (i.e. $\mathbb{T}^{\ast
}=\mathbb{R-}\left\{  0\right\}  $) since%
\[
\delta_{\pm}(s,0)=0\text{ for all }s\in\lbrack1,\infty)\text{.}%
\]
In the table below, we state different time scales with their corresponding
shift operators.
\[%
\begin{tabular}
[c]{|c||c|c|c|c|c|}\hline
$\mathbb{T}$ & $t_{0}$ & $t^{\ast}$ & $\mathbb{T}^{\ast}$ & $\delta_{-}(s,t)$
& $\delta_{+}(s,t)$\\\hline\hline
$\mathbb{R}$ & $0$ & N/A & $\mathbb{R}$ & $t-s$ & $t+s$\\\hline
$\mathbb{Z}$ & $0$ & N/A & $\mathbb{Z}$ & $t-s$ & $t+s$\\\hline
$q^{\mathbb{Z}}\cup\left\{  0\right\}  $ & $1$ & $0$ & $q^{\mathbb{Z}}$ &
$\frac{t}{s}$ & $st$\\\hline
$\mathbb{N}^{1/2}$ & $0$ & N/A & $\mathbb{N}^{1/2}$ & $\sqrt{t^{2}-s^{2}}$ &
$\sqrt{t^{2}+s^{2}}$\\\hline
\end{tabular}
\ \ \ \ \ \ \ \
\]

\end{example}

The proof of the next lemma is a direct consequence of Definition \ref{shift}.

\begin{lemma}
\label{lem pro} Let $\delta_{-}$ and $\delta_{+}$ be the shift operators
associated with the initial point $t_{0}$. We have

\begin{enumerate}
\item[i.] $\delta_{-}(t,t)=t_{0}$ for all $t\in\lbrack t_{0},\infty
)_{\mathbb{T}}.$

\item[ii.] $\delta_{-}(t_{0},t)=t$ for all $t\in\mathbb{T}^{\ast},$

\item[iii.] If $(s,t)\in$ $\mathcal{D}_{+}$, then $\delta_{+}(s,t)=u$ implies
$\delta_{-}(s,u)=t$. Conversely, if $(s,u)\in$ $\mathcal{D}_{-}$, then
$\delta_{-}(s,u)=t$ implies $\delta_{+}(s,t)=u$.

\item[iv.] $\delta_{+}(t,\delta_{-}(s,t_{0}))=\delta_{-}(s,t)$ for all
$(s,t)\in$ $\mathcal{D}(\delta_{+})$ with $t\geq t_{0,}$

\item[v.] $\delta_{+}(u,t)=\delta_{+}(t,u)$ for all $(u,t)\in\left(  \lbrack
t_{0},\infty)_{\mathbb{T}}\times\lbrack t_{0},\infty)_{\mathbb{T}}\right)
\cap\mathcal{D}_{+}$

\item[vi.] $\delta_{+}(s,t)\in\lbrack t_{0},\infty)_{\mathbb{T}}$ for all
$(s,t)\in$ $\mathcal{D}_{+}$ with $t\geq t_{0,}$,

\item[vii.] $\delta_{-}(s,t)\in\lbrack t_{0},\infty)_{\mathbb{T}}$ for all
$(s,t)\in$ $\left(  [t_{0},\infty)_{\mathbb{T}}\times\lbrack s,\infty
)_{\mathbb{T}}\right)  \cap\mathcal{D}_{-,}$

\item[viii.] If $\delta_{+}(s,.)$ is $\Delta-$differentiable in its second
variable, then $\delta_{+}^{\Delta_{t}}(s,.)>0$,

\item[ix.] $\delta_{+}(\delta_{-}(u,s),\delta_{-}(s,v))=\delta_{-}(u,v)$ for
all $(s,v)\in\left(  \lbrack t_{0},\infty)_{\mathbb{T}}\times\lbrack
s,\infty)_{\mathbb{T}}\right)  \cap\mathcal{D}_{-}$ and $(u,s)\in\left(
\lbrack t_{0},\infty)_{\mathbb{T}}\times\lbrack u,\infty)_{\mathbb{T}}\right)
\cap\mathcal{D}_{-}$,

\item[x.] If $(s,t)\in\mathcal{D}_{-}$ and $\delta_{-}(s,t)=t_{0}$, then $s=t$.
\end{enumerate}
\end{lemma}

\begin{proof}
(i) is obtained from P.3-5 since%
\[
\delta_{-}(t,t)=\delta_{-}(t,\delta_{+}(t,t_{0}))=t_{0}\text{ for all }%
t\in\mathbb{T}^{\ast}.
\]
(ii) is obtained from P.3-P.4 since%
\[
\delta_{-}(t_{0},t)=\delta_{-}(t_{0},\delta_{+}(t_{0},t))=t.
\]
Let $u:=\delta_{+}(s,t)$. By P.4 we have $(s,u)\in\mathcal{D}_{-}$ for all
$(s,t)\in\mathcal{D}_{+}$, and hence,%
\[
\delta_{-}(s,u)=\delta_{-}(s,\delta_{+}(s,t))=t.
\]
The latter part of (iii) can be done similarly.. We have (iv) since P.3 and
P.5 yield%
\[
\delta_{+}(t,\delta_{-}(s,t_{0}))=\delta_{-}(s,\delta_{+}(t,t_{0}))=\delta
_{-}(s,t).
\]
P.3 and P.5 guarantee that%
\[
t=\delta_{+}(t,t_{0})=\delta_{+}(t,\delta_{-}(u,u))=\delta_{-}(u,\delta
_{+}(t,u))
\]
for all $(u,t)\in\left(  \lbrack t_{0},\infty)_{\mathbb{T}}\times\lbrack
t_{0},\infty)_{\mathbb{T}}\right)  \cap\mathcal{D}_{+}$. Using (iii) we have%
\[
\delta_{+}(u,t)=\delta_{+}(u,\delta_{-}(u,\delta_{+}(t,u)))=\delta_{+}(t,u).
\]
This proves (v). To prove (vi) and (vii) we use P.1-2 to get%
\[
\delta_{+}(s,t)\geq\delta_{+}(t_{0},t)=t\geq t_{0}%
\]
for all $(s,t)\in$ $\left(  [t_{0},\infty)\times\lbrack t_{0},\infty
)_{\mathbb{T}}\right)  \cap\mathcal{D}_{+}$ and%
\[
\delta_{-}(s,t)\geq\delta_{-}(s,s)=t_{0}%
\]
for all $(s,t)\in\left(  \lbrack t_{0},\infty)_{\mathbb{T}}\times\lbrack
s,\infty)_{\mathbb{T}}\right)  \cap\mathcal{D}_{-}$. Since $\delta_{+}(s,t)$
is strictly increasing in its second variable we have (viii) by
\cite[Corollary 1.16]{book2}. (ix) is proven as follows: from P.5 and (v) we
have%
\begin{align*}
\delta_{+}(\delta_{-}(u,s),\delta_{-}(s,v))  &  =\delta_{-}(s,\delta
_{+}(v,\delta_{-}(u,s)))\\
&  =\delta_{-}(s,\delta_{-}(u,\delta_{+}(v,s)))\\
&  =\delta_{-}(s,\delta_{+}(s,\delta_{-}(u,v)))\\
&  =\delta_{-}(u,v)
\end{align*}
for all $(s,v)\in\left(  \lbrack t_{0},\infty)_{\mathbb{T}}\times\lbrack
s,\infty)_{\mathbb{T}}\right)  \cap\mathcal{D}_{-}$ and $(u,s)\in\left(
\lbrack t_{0},\infty)_{\mathbb{T}}\times\lbrack u,\infty)_{\mathbb{T}}\right)
\cap\mathcal{D}_{-}$. Suppose $(s,t)\in\mathcal{D}_{-}$ $=\left\{
(s,t)\in\lbrack t_{0},\infty)_{\mathbb{T}}\times\mathbb{T}^{\ast}:\delta
_{-}(s,t)\in\mathbb{T}^{\ast}\right\}  $ and $\delta_{-}(s,t)=t_{0}$. Then by
P.4 we have%
\[
t=\delta_{+}(s,\delta_{-}(s,t))\in\delta_{+}(s,t_{0})=s.
\]
This is (x). The proof is complete.
\end{proof}

Notice that the shift operators $\delta_{\pm}$ are defined once the initial
point $t_{0}\in\mathbb{T}^{\ast}$ is known. For instance, we choose the
initial point $t_{0}=0$ to define shift operators $\delta_{\pm}(s,t)=t\pm s$
on $\mathbb{T}=\mathbb{R}$. However, if we choose $\lambda\in(0,\infty)$ as
the initial point, then the new shift operators associated with $\lambda$ are
defined by $\widetilde{\delta}_{\pm}(s,t)=t\mp\lambda\pm s$. In terms of
$\delta_{\pm}$ the operators $\widetilde{\delta}_{\pm}$ can be given as%
\[
\widetilde{\delta}_{\pm}(s,t)=\delta_{\mp}(\lambda,\delta_{\pm}(s,t)).
\]

\begin{example}
In the following, we give some particular time scales to show the change in
the formula of shift operators as the initial point changes.%
\[%
\begin{tabular}
[c]{c||cc|cc|cc}
& \multicolumn{2}{||c|}{$\mathbb{T}=\mathbb{N}^{1/2}$} &
\multicolumn{2}{|c|}{$\mathbb{T}=h\mathbb{Z}$} &
\multicolumn{2}{|c}{$\mathbb{T}=2^{\mathbb{N}}$}\\\hline\hline
$t_{0}$ & $0$ & $\lambda$ & $0$ & $h\lambda$ & $1$ & $2^{\lambda}$\\
$\delta_{-}(s,t)$ & $\sqrt{t^{2}-s^{2}}$ & $\sqrt{t^{2}+\lambda^{2}-s^{2}}$ &
$t-s$ & $t+h\lambda-s$ & $t/s$ & $2^{\lambda}ts^{-1}$\\
$\delta_{+}(s,t)$ & $\sqrt{t^{2}+s^{2}}$ & $\sqrt{t^{2}-\lambda^{2}+s^{2}}$ &
$t+s$ & $t-h\lambda+s$ & $ts$ & $2^{-\lambda}ts$%
\end{tabular}
\]
where $\lambda\in\mathbb{Z}_{+}$, $\mathbb{N}^{1/2}:=\{\sqrt{n}:n\in
\mathbb{N\}}$, $2^{\mathbb{N}}:=\{2^{n}:n\in\mathbb{N\}}$, and $h\mathbb{Z}%
:\mathbb{=\{}hn:$ $n\in\mathbb{Z\}}$.
\end{example}

\section{Delay function}

In this section we introduce the delay function on time scales that will be
used for the construction of the Lyapunov functional.

\begin{definition}
Let $\mathbb{T}$ be a time scale that is unbounded above\ and $t_{0}%
\in\mathbb{T}^{\ast}$ an element such that there exist the shift operators
$\delta_{\pm}:[t_{0},\infty)\times\mathbb{T}^{\ast}\rightarrow\mathbb{T}%
^{\ast}$ associated with $t_{0}$. Suppose that $h\in(t_{0},\infty
)_{\mathbb{T}}$ is a constant such that $(h,t)\in D_{\pm}$ for all
$t\in\lbrack t_{0},\infty)_{\mathbb{T}}$, the function $\delta_{-}(h,t)$ is
differentiable with an $rd$-continuous derivative, and $\delta_{-}(h,t)$ maps
$[t_{0},\infty)_{\mathbb{T}}$ onto $[\delta_{-}(h,t_{0}),\infty)_{\mathbb{T}}%
$. Then the function $\delta_{-}(h,t)$ is called the delay function generated
by the shift $\delta_{-}$ on the time scale $\mathbb{T}$.
\end{definition}

It is obvious from P.2 and (iii) of Lemma \ref{lem pro} that%
\begin{equation}
\delta_{-}(h,t)<\delta_{-}(t_{0},t)=t\text{ for all }t\in\lbrack t_{0}%
,\infty)_{\mathbb{T}}\text{.} \label{delay less}%
\end{equation}
Notice that $\delta_{-}(h,.)$ is strictly increasing and it is invertible.
Hence, by P.4-5 $\delta_{-}^{-1}(h,t)=\delta_{+}(h,t)$.

Hereafter, we shall suppose that $\mathbb{T}$ is a time scale with the delay
function $\delta_{-}(h,.):[t_{0},\infty)_{\mathbb{T}}\rightarrow\lbrack
\delta_{-}(h,t_{0}),\infty)_{\mathbb{T}}$, where $t_{0}\in\mathbb{T}$ is
fixed. Denote by $\mathbb{T}_{1}$ and $\mathbb{T}_{2}$ the sets%
\begin{equation}
\mathbb{T}_{1}=[t_{0},\infty)_{\mathbb{T}}\text{\ \ and }\mathbb{T}_{2}%
=\delta_{-}(h,\mathbb{T}_{1}).\label{T12}%
\end{equation}
Evidently, $\mathbb{T}_{1}$ is closed in $\mathbb{R}$. By definition we have
$\mathbb{T}_{2}=[\delta_{-}(h,t_{0}),\infty)_{\mathbb{T}}$. Hence,
$\mathbb{T}_{1}$ and $\mathbb{T}_{2}$ are both time scales. Let $\sigma_{1}$
and $\sigma_{2}$ denote the forward jump operators on the time scales
$\mathbb{T}_{1}$ and $\mathbb{T}_{2}$, respectively. By (\ref{delay less}%
-\ref{T12})
\[
\mathbb{T}_{1}\subset\mathbb{T}_{2}\subset\mathbb{T}.
\]
Thus,%
\[
\sigma(t)=\sigma_{2}(t)\text{ for all }t\in\mathbb{T}_{2}%
\]
and%
\[
\sigma(t)=\sigma_{1}(t)=\sigma_{2}(t)\text{ for all }t\in\mathbb{T}_{1}.
\]
That is, $\sigma_{1}$ and $\sigma_{2}$ are the restrictions of the forward
jump operator $\sigma:\mathbb{T\rightarrow T}$ to the time scales
$\mathbb{T}_{1}$ and $\mathbb{T}_{2}$, respectively, i.e.,%
\[
\sigma_{1}=\left.  \sigma\right\vert _{\mathbb{T}_{1}}\text{ and }\sigma
_{2}=\left.  \sigma\right\vert _{\mathbb{T}_{2}}\text{.}%
\]
Recall that the Hilger derivatives $\Delta$, $\Delta_{1}$, and $\Delta_{2}$ on
the time scales $\mathbb{T}$, $\mathbb{T}_{1}$, and $\mathbb{T}_{2}$ are
defined in terms of the forward jumps $\sigma$, $\sigma_{1}$, and $\sigma_{2}%
$, respectively. Hence, if $f$ is a differentiable function at $t\in
\mathbb{T}_{2}$, then we have%
\[
f^{\Delta_{2}}(t)=f^{\Delta_{1}}(t)=f^{\Delta}(t),\ \ \text{\ for all }%
t\in\mathbb{T}_{1}.
\]
Similarly, if $a,b\in\mathbb{T}_{2}$ are two points with $a<b$ and if $f$ is a
$rd$-continuous function on the interval $(a,b)_{\mathbb{T}_{2}}$, then%
\[
\int_{a}^{b}f(s)\Delta_{2}s=\int_{a}^{b}f(s)\Delta s.
\]
The next result is essential for future calculations.

\begin{lemma}
The delay function $\delta_{-}(h,t)$ preserves the structure of the points in
$\mathbb{T}_{1}$. That is,%
\[
\sigma_{1}(\widehat{t})=\widehat{t}\text{ implies }\sigma_{2}(\delta
_{-}(h,\widehat{t}))=\delta_{-}(h,\widehat{t}).
\]%
\[
\sigma_{1}(\widehat{t})>\widehat{t}\text{ implies }\sigma_{2}(\delta
_{-}(h,\widehat{t})>\delta_{-}(h,\widehat{t}).
\]

\end{lemma}

\begin{proof}
By definition $\sigma_{1}(t)\geq t$ for all $t\in\mathbb{T}_{1}$. Thus,%
\[
\delta_{-}(h,\sigma_{1}(t))\geq\delta_{-}(h,t).
\]
Since $\sigma_{2}(\delta_{-}(h,t))$ is the smallest element satisfying%
\[
\sigma_{2}(\delta_{-}(h,t))\geq\delta_{-}(h,t),
\]
we get%
\begin{equation}
\delta_{-}(h,\sigma_{1}(t))\geq\sigma_{2}(\delta_{-}(h,t))\text{ for all }%
t\in\mathbb{T}_{1}\text{.} \label{1}%
\end{equation}
If $\sigma_{1}(\widehat{t})=\widehat{t}$, then we have%
\[
\delta_{-}(h,\widehat{t})=\delta_{-}(h,\sigma_{1}(\widehat{t}))\geq\sigma
_{2}(\delta_{-}(h,\widehat{t})).
\]
That is,%
\[
\delta_{-}(h,\widehat{t})=\sigma_{2}(\delta_{-}(h,\widehat{t})).
\]
If $\sigma_{1}(\widehat{t})>\widehat{t}$, then%
\[
(\widehat{t},\sigma_{1}(\widehat{t}))_{\mathbb{T}_{1}}=(\widehat{t},\sigma
_{1}(\widehat{t}))_{\mathbb{T}}=\varnothing
\]
and%
\[
\delta_{-}(h,\sigma_{1}(\widehat{t}))>\delta_{-}(h,\widehat{t}).
\]
Suppose the contrary. That is $\delta_{-}(h,\widehat{t})$ is right dense;
namely $\sigma_{2}(\delta_{-}(h,\widehat{t}))=\delta_{-}(h,\widehat{t})$. This
along with (\ref{1}) implies%
\[
(\delta_{-}(h,\widehat{t}),\delta_{-}(h,\sigma_{1}(\widehat{t})))_{\mathbb{T}%
_{2}}\neq\varnothing\text{.}%
\]
Pick one element $s\in(\delta_{-}(h,\widehat{t}),\delta_{-}(h,\sigma
_{1}(\widehat{t})))_{\mathbb{T}_{2}}$. Since $\delta_{-}(h,t)$ is strictly
increasing in $t$ and invertible, there should be an element $t\in(\widehat
{t},\sigma_{1}(\widehat{t}))_{\mathbb{T}_{1}}$ such that $\delta_{-}(h,t)=s$.
This leads to a contradiction. Hence, $\delta_{-}(h,\widehat{t})$ must be
right scattered.
\end{proof}

Using the preceding lemma and applying the fact that $\sigma_{2}(u)=\sigma(u)$
for all $u\in\mathbb{T}_{2}$ we arrive at the following result.

\begin{corollary}
\label{Cor 1} We have%
\[
\delta_{-}(h,\sigma_{1}(t))=\sigma_{2}(\delta_{-}(h,t))\text{ for all }%
t\in\mathbb{T}_{1}\text{.}%
\]
Thus,%
\begin{equation}
\delta_{-}(h,\sigma(t))=\sigma(\delta_{-}(h,t))\text{ for all }t\in
\mathbb{T}_{1}\text{.} \label{sigma delta}%
\end{equation}

\end{corollary}

By (\ref{sigma delta}) we have%
\[
\delta_{-}(h,\sigma(s))=\sigma(\delta_{-}(h,s))\text{ for all }s\in\lbrack
t_{0},\infty)_{\mathbb{T}}\text{.}%
\]
Substituting $s=\delta_{+}(h,t)$ we obtain%
\[
\delta_{-}(h,\sigma(\delta_{+}(h,t)))=\sigma(\delta_{-}(h,\delta
_{+}(h,t)))=\sigma(t)\text{.}%
\]
This and (iv) of Lemma \ref{lem pro} imply%
\[
\sigma(\delta_{+}(h,t))=\delta_{+}(h,\sigma(t))\text{ for all }t\in
\lbrack\delta_{-}(h,t_{0}),\infty)_{\mathbb{T}}\text{.}%
\]

\begin{example}
In the following, we give some time scales with their shift operators:
\[%
\begin{tabular}
[c]{|c||c|c|c|}\hline
$\mathbb{T}$ & $h$ & $\delta_{-}(h,t)$ & $\delta_{+}(h,t)$\\\hline\hline
$\mathbb{R}$ & $\in\mathbb{R}_{+}$ & $t-h$ & $t+h$\\\hline
$\mathbb{Z}$ & $\in\mathbb{Z}_{+}$ & $t-h$ & $t+h$\\\hline
$q^{\mathbb{Z}}\cup\left\{  0\right\}  $ & $\in q^{\mathbb{Z}_{+}}$ &
$\frac{t}{h}$ & $ht$\\\hline
$\mathbb{N}^{1/2}$ & $\in\mathbb{Z}_{+}$ & $\sqrt{t^{2}-h^{2}}$ & $\sqrt
{t^{2}+h^{2}}$\\\hline
\end{tabular}
\]

\end{example}

\begin{example}
There is no delay function $\delta_{-}(h,.):[0,\infty)_{\widetilde{\mathbb{T}%
}}\rightarrow\lbrack\delta_{-}(h,0),\infty)_{\mathbb{T}}$ on the time scale
$\widetilde{\mathbb{T}}\mathbb{=(-\infty},0]\cup\lbrack1,\infty)$.\newline
Suppose the contrary that there exists such a delay function on $\widetilde
{\mathbb{T}}$. Then since $0$ is right scattered in $\widetilde{\mathbb{T}%
}_{1}:=[0,\infty)_{\widetilde{\mathbb{T}}}$ the point $\delta_{-}(h,0)$ must
be right scattered in $\widetilde{\mathbb{T}}_{2}=[\delta_{-}(h,0),\infty
)_{\mathbb{T}}$, i.e., $\sigma_{2}(\delta_{-}(h,0))>\delta_{-}(h,0)$. Since
$\sigma_{2}(t)=\sigma(t)$ for all $t\in\lbrack\delta_{-}(h,0),0)_{\mathbb{T}}%
$, we have%
\[
\sigma(\delta_{-}(h,0))=\sigma_{2}(\delta_{-}(h,0))>\delta_{-}(h,0).
\]
That is, $\delta_{-}(h,0)$ must be right scattered in $\widetilde{\mathbb{T}}%
$. However, in $\widetilde{\mathbb{T}}$ we have $\delta_{-}(h,0)<0$, that is,
$\delta_{-}(h,0)$ is right dense. This leads to a contradiction.
\end{example}

\begin{theorem}
(\textbf{Substitution}) \label{thm2.2} \cite[Theorem 1.98]{book} Assume
$\nu:\mathbb{T}\rightarrow\mathbb{R}$ is strictly increasing and
$\tilde{\mathbb{T}}:=\nu(\mathbb{T})$ is a time scale. If $f:\mathbb{T}%
\rightarrow\mathbb{R}$ is an rd-continuous function and $\nu$ is
differentiable with rd-continuous derivative, then for $a,b\in\mathbb{T}$,
\begin{equation}
\int_{a}^{b}\!g(t,s)\nu^{\Delta}(s)\,\Delta s=\int_{\nu(a)}^{\nu(b)}%
g(t,\nu^{-1}(s))\,\tilde{\Delta}s. \label{substitute}%
\end{equation}

\end{theorem}

First, since the operator $\delta:[t_{0},\infty)_{\mathbb{T}}\rightarrow$
$[\delta(t_{0}),\infty)_{\mathbb{T}}$ is strictly increasing, it is bijection.
If we substitute $\nu(t)=\delta_{-}(h,t)$ and%
\[
f(t,s)=g(t,\delta_{-}^{-1}(h,s))=g(t,\delta_{+}(h,s))
\]
into (\ref{substitute}), we obtain%
\begin{equation}
\int_{a}^{b}\!f(t,\delta_{-}(h,s))\delta_{-}^{\Delta_{1}}(h,s)\,\Delta
_{1}s=\int_{\delta_{-}(h,a)}^{\delta_{-}(h,b)}f(t,s)\,\Delta_{2}s \label{5}%
\end{equation}
for $a,b\in\mathbb{T}_{1}$. For any $t\in\mathbb{T}_{1}$, we have $[\delta
_{-}(h,t_{0}),t)_{\mathbb{T}_{1}}\subset\mathbb{T}_{2}$. This and (\ref{5})
yield
\begin{align}
\int_{\delta_{-}(h,t)}^{t}f(t,s)\Delta s  &  =\int_{\delta_{-}(h,t)}%
^{t}f(t,s)\Delta_{2}s\nonumber\\
&  =\int_{\delta_{-}(h,t)}^{\delta_{-}(h,t_{0})}f(t,s)\Delta_{2}s+\int
_{\delta_{-}(h,t_{0})}^{t}f(t,s)\Delta_{2}s\nonumber\\
&  =\int_{t}^{t_{0}}f(t,\delta_{-}(h,s))\delta_{-}^{\Delta_{1}}(h,s)\,\Delta
_{1}s+\int_{\delta_{-}(h,t_{0})}^{t}f(t,s)\Delta s\nonumber\\
&  =\int_{t}^{t_{0}}f(t,\delta_{-}(h,s))\delta_{-}^{\Delta}(h,s)\,\Delta
s+\int_{\delta_{-}(h,t_{0})}^{t}f(t,s)\Delta s. \label{3}%
\end{align}
The the formula%
\begin{align}
\left[  \int_{\delta_{-}(h,t)}^{t}f(t,s)\Delta s\right]  ^{\Delta}  &
=f(\sigma(t),t)-f(\sigma(t),\delta_{-}(h,t))\delta_{-}^{\Delta}%
(h,t)\nonumber\\
&  +\int_{\delta_{-}(h,t)}^{t}f^{\Delta}(t,s)\Delta s \label{4}%
\end{align}
follows from (\ref{3}) and Theorem \ref{theorem 1.117}.

\begin{theorem}
\label{order of integration} Let $k$ be an $rd$-continuous function. Then%
\begin{equation}
\int_{\delta_{-}(h,t)}^{t}\Delta s\int_{s}^{t}k(u)\Delta u=\int_{\delta
_{-}(h,t)}^{t}\Delta u\int_{\delta_{-}(h,t)}^{\sigma(u)}k(u)\Delta s.\label{8}%
\end{equation}

\end{theorem}

\begin{proof}
Substituting%
\[
f(s)=s-\delta_{-}(h,t),\ \ \ \ g(s)=\int_{s}^{t}k(u)\Delta u
\]
into the formula%
\[
\int_{a}^{z}f(\sigma(x))g(x)\Delta x=\left[  f(x)g(x)\right]  _{a}^{z}%
-\int_{a}^{z}f^{\Delta}(x)g(x)\Delta x
\]
(see \cite[Theorem 1.77]{book}) and using Lemma \ref{lem pro} we get%
\begin{align}
\int_{\delta_{-}(h,t)}^{t}\Delta s\int_{s}^{t}k(u)\Delta u &  =\int
_{\delta_{-}(h,t)}^{t}\left[  \sigma(s)-\delta_{-}(h,t)\right]  k(s)\Delta
s\nonumber\\
&  =\int_{\delta_{-}(h,t)}^{t}\Delta u\int_{\delta_{-}(h,t)}^{\sigma
(u)}k(u)\Delta s.\label{8a}%
\end{align}

\end{proof}

\section{Stability analysis using Lyapunov's method}

Let $\mathbb{T}$ be a time scale having a delay function $\delta_{-}(h,t)$
where $h\geq t_{0}$ and $t_{0}\in\mathbb{T}$ is nonnegative and fixed. In this
section we consider the equation%
\begin{equation}
x^{\Delta}(t)=a(t)x(t)+b(t)x(\delta_{-}(h,t))\delta_{-}^{\Delta}%
(h,t),\ \ \ t\in\lbrack t_{0},\infty)_{\mathbb{T}} \label{9}%
\end{equation}
and assume that%
\begin{equation}
\left\vert \delta_{-}^{\Delta}(h,t)\right\vert \leq M<\infty\ \text{ for all
}t\in\lbrack t_{0},\infty)_{\mathbb{T}}. \label{delta turev}%
\end{equation}
Let $\psi$: $[\delta_{-}(h,t_{0}),t_{0}]_{\mathbb{T}}\rightarrow\mathbb{R}$ be
$rd$-continuous and let $x(t):=x(t,t_{0},\psi)$ be the solution of Eq.
(\ref{9}) on $[t_{0},\infty)_{\mathbb{T}}$ with $x(t)=\psi(t)$ on $[\delta
_{-}(h,t_{0}),t_{0}]_{\mathbb{T}}$. Let $\left\Vert \varphi\right\Vert
=\sup\left\{  |\varphi(t)|:t\in\lbrack\delta_{-}(h,t_{0}),t_{0})_{\mathbb{T}%
}\right\}  $.

Observe that using (\ref{4}) Eq. (\ref{9}) can be rewritten as follows%
\begin{equation}
x^{\Delta}(t)=Q(t)x(t)-\left[  \int_{\delta_{-}(h,t)}^{t}b(\delta
_{+}(h,s))x(s)\Delta s\right]  ^{\Delta_{t}}, \label{10}%
\end{equation}
where
\[
Q(t):=a(t)+b(\delta_{+}(h,t))
\]
and $\Delta_{t}$ indicates the delta derivative with respect to $t$.

\begin{lemma}
\label{lem2} Let%
\begin{equation}
A(t):=x(t)+%
{\displaystyle\int\limits_{\delta_{-}(h,t)}^{t}}
b(\delta_{+}(h,s))x(s)\Delta s\label{11a}%
\end{equation}
and
\begin{equation}
\beta(t):=t-\delta_{-}(h,t).\label{11b}%
\end{equation}
Assume that there exists a $\lambda>0$ such that%
\begin{equation}
-\frac{\lambda\delta_{-}^{\Delta}(h,t)}{\beta(t)+\lambda\left[  \beta
(t)+\mu(t)\right]  }\leq Q(t)\leq-\lambda\left[  \beta(t)+\mu(t)\right]
b(\delta_{+}(h,t))^{2}-\mu(t)Q^{2}(t)\label{11}%
\end{equation}
for all $t\in\lbrack t_{0},\infty)_{\mathbb{T}}$. If%
\begin{equation}
V(t)=A(t)^{2}+\lambda%
{\displaystyle\int\limits_{\delta_{-}(h,t)}^{t}}
\Delta s%
{\displaystyle\int\limits_{s}^{t}}
b(\delta_{+}(h,u))^{2}x(u)^{2}\Delta u\label{12}%
\end{equation}
then, along the solutions of Eq. (\ref{9}) we have%
\begin{equation}
V^{\Delta}(t)\leq Q(t)V(t)\text{ for all }t\in\lbrack t_{0},\infty
)_{\mathbb{T}}\text{.}\label{13}%
\end{equation}

\end{lemma}

\begin{proof}
It is obvious from (\ref{10}) and (\ref{11a}) that%
\[
A^{\Delta}(t)=Q(t)x(t).
\]
Then by (\ref{4}) and the formula $A(\sigma(t))=A(t)+\mu(t)A(t)$ we have%
\begin{align*}
V^{\Delta}(t) &  =\left[  A(t)+A(\sigma(t))\right]  A^{\Delta}(t)+\lambda
\int_{t}^{\sigma(t)}b(\delta_{+}(h,u))^{2}x(u)^{2}\Delta u\\
&  -\lambda\delta_{-}^{\Delta}(h,t)\int_{\delta_{-}(h,t)}^{\sigma(t)}%
b(\delta_{+}(h,u))^{2}x(u)^{2}\Delta u+\lambda\left(  t-\delta_{-}%
(h,t)\right)  b(\delta_{+}(h,t))^{2}x(t)^{2}\\
&  =\left[  2A(t)+\mu(t)Q(t)x(t)\right]  Q(t)x(t)-\lambda\delta_{-}^{\Delta
}(h,t)\int_{\delta_{-}(h,t)}^{\sigma(t)}b(\delta_{+}(h,u))^{2}x(u)^{2}\Delta
u\\
&  +\lambda\left[  \beta(t)+\mu(t)\right]  b(\delta_{+}(h,t))^{2}x(t)^{2}.
\end{align*}
Using the identity%
\begin{equation}
2A(t)x(t)=x^{2}(t)+A^{2}(t)-\left(
{\displaystyle\int\limits_{\delta_{-}(h,t)}^{t}}
b(\delta_{+}(h,s))x(s)\Delta s\right)  ^{2}\label{A_Q}%
\end{equation}
and condition (\ref{11}) we have%
\begin{align}
V^{\Delta}(t) &  =Q(t)V(t)+R(t)\nonumber\\
&  +x^{2}(t)\left[  \lambda\left(  \beta(t)+\mu(t)\right)  b(\delta
_{+}(h,t))^{2}+Q(t)+\mu(t)Q^{2}(t)\right]  \nonumber\\
&  \leq Q(t)V(t)+R(t),\label{v(t)}%
\end{align}
where%
\begin{align}
R(t) &  =-Q(t)\left(
{\displaystyle\int\limits_{\delta_{-}(h,t)}^{t}}
b(\delta_{+}(h,s))x(s)\Delta s\right)  ^{2}\nonumber\\
&  -\lambda\delta_{-}^{\Delta}(h,t)\int_{\delta_{-}(h,t)}^{\sigma(t)}%
b(\delta_{+}(h,u))^{2}x(u)^{2}\Delta u\nonumber\\
&  -\lambda Q(t)%
{\displaystyle\int\limits_{\delta_{-}(h,t)}^{t}}
\Delta s%
{\displaystyle\int\limits_{s}^{t}}
b(\delta_{+}(h,u))^{2}x(u)^{2}\Delta u.\label{r(t)}%
\end{align}
Hereafter, we will show that (\ref{11}) implies $R(t)\leq0$. This and
(\ref{v(t)}) will enable us to derive the desired inequality (\ref{13}). First
we have
\begin{equation}
\int_{\delta_{-}(h,t)}^{\sigma(t)}b(\delta_{+}(h,u))^{2}x(u)^{2}\Delta
u\geq\int_{\delta_{-}(h,t)}^{t}b(\delta_{+}(h,u))^{2}x(u)^{2}\Delta
u.\label{r1}%
\end{equation}
From H\"{o}lder's inequality \cite[Theorem 6.13]{book} we get%
\begin{equation}
\left(
{\displaystyle\int\limits_{\delta_{-}(h,t)}^{t}}
b(\delta_{+}(h,s))x(s)\Delta s\right)  ^{2}\leq\beta(t)%
{\displaystyle\int\limits_{\delta_{-}(h,t)}^{t}}
b(\delta_{+}(h,s))^{2}x(s)^{2}\Delta s.\label{holder}%
\end{equation}
On the other hand, (\ref{8}) yields
\begin{align}%
{\displaystyle\int\limits_{\delta_{-}(h,t)}^{t}}
\Delta s%
{\displaystyle\int\limits_{s}^{t}}
b(\delta_{+}(h,u))^{2}x(u)^{2}\Delta u &  =%
{\displaystyle\int\limits_{\delta_{-}(h,t)}^{t}}
\Delta u%
{\displaystyle\int\limits_{\delta_{-}(h,t)}^{\sigma(u)}}
b(\delta_{+}(h,u))^{2}x(u)^{2}\Delta s\nonumber\\
&  =%
{\displaystyle\int\limits_{\delta_{-}(h,t)}^{t}}
\left[  \sigma(u)-\delta_{-}(h,t)\right]  b(\delta_{+}(h,u))^{2}x(u)^{2}\Delta
u\nonumber\\
&  \leq\left[  \beta(t)+\mu(t)\right]
{\displaystyle\int\limits_{\delta_{-}(h,t)}^{t}}
b(\delta_{+}(h,u))^{2}x(u)^{2}\Delta u.\label{holder1}%
\end{align}
Substituting (\ref{holder}) and (\ref{holder1}) into (\ref{r(t)}) and using
(\ref{r1}) together with $\delta_{-}^{\Delta}(h,t)>0$ we deduce%
\[
R(t)\leq-\left\{  \left(  \beta(t)+\left[  \lambda\beta(t)+\mu(t)\right]
\right)  Q(t)+\lambda\delta_{-}^{\Delta}(h,t)\right\}
{\displaystyle\int\limits_{\delta_{-}(h,t)}^{t}}
b(\delta_{+}(h,s))^{2}x(s)^{2}\Delta s.
\]
Hence, using the left-hand side of (\ref{11}) we arrive at the inequality
$R(t)\leq0$. The proof is complete.
\end{proof}

In preparation for the proof of the next theorem we state the following lemma.

\begin{lemma}
\label{Lemma ep}If $\varphi\in\mathcal{R}^{+}$, then%
\begin{equation}
0<e_{\varphi}(t,s)\leq\exp\left(  \int_{s}^{t}\varphi(r)\Delta r\right)
\label{zero exp}%
\end{equation}
for all $t\in\lbrack s,\infty)_{\mathbb{T}}$.
\end{lemma}

\begin{theorem}
\label{thm2.1} Let $a\in\mathcal{R}^{+}$ and $Q\in\mathcal{R}$. Suppose the
hypothesis of Lemma \ref{lem2}. If there exists an $\alpha\in(t_{0}%
,h)_{\mathbb{T}}$ such that
\begin{equation}
(a,t)\in\mathcal{D}_{\pm}\ \ \ \text{for all\ \ \ }t\in\lbrack t_{0}%
,\infty)_{\mathbb{T}}\label{2.8a}%
\end{equation}
and%
\begin{equation}
\delta_{-}(h,t)\leq\frac{\delta_{-}(\alpha,t)+\delta_{-}(h,\delta_{-}%
(\alpha,t))}{2}\text{ for all }t\in\lbrack\alpha,\infty)_{\mathbb{T}%
},\label{2.8}%
\end{equation}
then any solution $x(t)=x(t,t_{0},\varphi)$ of (\ref{9}) satisfies the
exponential inequalities%
\begin{equation}
\left\vert x(t)\right\vert \leq\sqrt{\frac{2}{\left(  1-\frac{1}{\xi
(t)}\right)  }V(t_{0})}e^{\frac{1}{2}\int_{t_{0}}^{\delta_{-}(\alpha
,t)}Q(s)\Delta s}\label{2.9}%
\end{equation}
for all $t\in\lbrack\alpha,\infty)_{\mathbb{T}}$ and
\begin{equation}
\left\vert x(t)\right\vert \leq\left\Vert \psi\right\Vert e^{\int_{t_{0}}%
^{t}a(s)\Delta s}\left[  1+M\int_{t_{0}}^{t}\left\vert \frac{b(s)}%
{1+\mu(s)a(s)}\right\vert e^{-\int_{t_{0}}^{s}a(u)\Delta u}\Delta s\right]
\label{2.10}%
\end{equation}
for all $t\in\lbrack t_{0},\alpha)_{\mathbb{T}}$, where $M$ is as defined by
(\ref{delta turev}),
\[
\xi(t):=1+\frac{\lambda\Lambda(t)}{\beta(t)}>1,
\]
and $\Lambda(t):=\delta_{-}(h,t)-\delta_{-}(h,\delta_{-}(\alpha,t)).$
\end{theorem}

\begin{proof}
Since $t_{0}<\alpha<h$ the condition (\ref{2.8}) implies%
\begin{equation}
\delta_{-}(h,t)<\delta_{-}(\alpha,t)\text{ for all }t\in\lbrack\alpha
,\infty)_{\mathbb{T}}\label{2.11}%
\end{equation}
and%
\begin{equation}
0<\Lambda(t)\leq\delta_{-}(\alpha,t)-\delta_{-}(h,t)\text{ for all }%
t\in\lbrack\alpha,\infty)_{\mathbb{T}},\label{2.12}%
\end{equation}
Let $V(t)$ be defined by (\ref{12}). First we get by (\ref{8a}), (\ref{12})
and (\ref{2.11}-\ref{2.12}) that%
\begin{align}
V(t) &  \geq\lambda%
{\displaystyle\int\limits_{\delta_{-}(h,t)}^{t}}
\Delta s%
{\displaystyle\int\limits_{s}^{t}}
b(\delta_{+}(h,u))^{2}x(u)^{2}\Delta u\nonumber\\
&  =\lambda\int_{\delta_{-}(h,t)}^{t}\left[  \sigma(u)-\delta_{-}(h,t)\right]
b(\delta_{+}(h,u))^{2}x(u)^{2}\Delta u\nonumber\\
&  \geq\lambda\int_{\delta_{-}(\alpha,t)}^{t}\left[  \sigma(u)-\delta
_{-}(h,t)\right]  b(\delta_{+}(h,u))^{2}x(u)^{2}\Delta u\nonumber\\
&  \geq\lambda\left[  \delta_{-}(\alpha,t)-\delta_{-}(h,t)\right]
\int_{\delta_{-}(\alpha,t)}^{t}b(\delta_{+}(h,u))^{2}x(u)^{2}\Delta
u.\nonumber
\end{align}
This along with (\ref{2.12}) yields
\begin{equation}
V(t)\geq\lambda\Lambda(t)\int_{\delta_{-}(\alpha,t)}^{t}b(\delta_{+}%
(h,u))^{2}x(u)^{2}\Delta u\label{v-1}%
\end{equation}
for all $t\in\lbrack\alpha,\infty)_{\mathbb{T}}$. Similarly, we get%
\begin{align}
V(\delta_{-}(\alpha,t)) &  \geq\lambda\int_{\delta_{-}(h,\delta_{-}%
(\alpha,t))}^{\delta_{-}(\alpha,t)}\left[  \sigma(u)-\delta_{-}(h,\delta
_{-}(\alpha,t))\right]  b(\delta_{+}(h,u))^{2}x(u)^{2}\Delta u\nonumber\\
&  \geq\lambda\int_{\delta_{-}(h,t)}^{\delta_{-}(\alpha,t)}\left[
\sigma(u)-\delta_{-}(h,\delta_{-}(\alpha,t))\right]  b(\delta_{+}%
(h,u))^{2}x(u)^{2}\Delta u\nonumber\\
&  \geq\lambda\Lambda(t)\int_{\delta_{-}(h,t)}^{\delta_{-}(\alpha,t)}%
b(\delta_{+}(h,u))^{2}x(u)^{2}\Delta u\label{v2}%
\end{align}
for all $t\in\lbrack\alpha,\infty)_{\mathbb{T}}$ since $\delta_{-}%
(\alpha,t)\leq\delta_{-}(t_{0},t)=t$. Utilizing (\ref{12}), (\ref{v-1}), and
(\ref{v2}) we obtain%
\begin{align}
V(t)+V(\delta_{-}(\alpha,t)) &  \geq A(t)^{2}+\lambda\Lambda(t)\int
_{\delta_{-}(\alpha,t)}^{t}b(\delta_{+}(h,u))^{2}x(u)^{2}\Delta u\nonumber\\
&  +\lambda\Lambda(t)\int_{\delta_{-}(h,t)}^{\delta_{-}(\alpha,t)}b(\delta
_{+}(h,u))^{2}x(u)^{2}\Delta u\nonumber\\
&  \geq A(t)^{2}+\lambda\Lambda(t)\int_{\delta_{-}(h,t)}^{t}b(\delta
_{+}(h,u))^{2}x(u)^{2}\Delta u\label{v3}%
\end{align}
for all $t\in\lbrack\alpha,\infty)_{\mathbb{T}}$. Substituting (\ref{holder})
and (\ref{11a}) into (\ref{v3}) we find%
\begin{align}
V(t)+V(\delta_{-}(\alpha,t)) &  \geq\left(  1-\frac{1}{\xi(t)}\right)
x^{2}(t)\nonumber\\
&  +\left[  \frac{1}{\sqrt{\xi(t)}}x(t)+\sqrt{\xi(t)}\left(
{\displaystyle\int\limits_{\delta_{-}(h,t)}^{t}}
b(\delta_{+}(h,u))x(u)\Delta u\right)  \right]  ^{2}\nonumber\\
&  \geq\left(  1-\frac{1}{\xi(t)}\right)  x^{2}(t)\label{v4}%
\end{align}
for all $t\in\lbrack\alpha,\infty)_{\mathbb{T}}$. Since $V^{\Delta}(t)\leq0$,
we get by (\ref{v4}) that%
\begin{equation}
\left(  1-\frac{1}{\xi(t)}\right)  x^{2}(t)\leq V(t)+V(\delta_{-}%
(\alpha,t))\leq2V(\delta_{-}(\alpha,t))\label{v5}%
\end{equation}
for all $t\in\lbrack\alpha,\infty)_{\mathbb{T}}$. Multiplying (\ref{13}) by
$e_{\ominus Q}(\sigma(s),t_{0})$ and integrating the resulting inequality from
$t_{0}$ to $t$ we derive%
\begin{align}
0 &  \geq\int_{t_{0}}^{t}\left[  V^{\Delta}(s)-Q(s)V(s)\right]  e_{\ominus
Q}(\sigma(s),t_{0})\Delta s\nonumber\\
&  =\int_{t_{0}}^{t}\left[  V(s)e_{\ominus Q}(s,t_{0})\right]  ^{\Delta}\Delta
s\nonumber\\
&  =V(t)e_{\ominus Q}(t,t_{0})-V(t_{0}).\label{vin}%
\end{align}
That is,%
\begin{equation}
V(t)\leq V(t_{0})e_{Q}(t,t_{0})\text{ for all }t\in\lbrack t_{0}%
,\infty)_{\mathbb{T}}\text{.}\label{v6}%
\end{equation}
Combining (\ref{v5}) and (\ref{v6}) we arrive at%
\[
x^{2}(t)\leq\frac{2}{\left(  1-\frac{1}{\xi(t)}\right)  }V(t_{0})e_{Q}%
(\delta_{-}(\alpha,t),t_{0})
\]
for all $t\in\lbrack\alpha,\infty)_{\mathbb{T}}$. The hypothesis
$Q\in\mathcal{R}$ and the condition (\ref{11}) guarantee that $Q(t)\in
\mathcal{R}^{+}$. Thus, (\ref{zero exp}) implies
\[
\left\vert x(t)\right\vert \leq\sqrt{\frac{2}{\left(  1-\frac{1}{\xi
(t)}\right)  }V(t_{0})}e^{\frac{1}{2}\int_{t_{0}}^{\delta_{-}(\alpha
,t)}Q(s)\Delta s}%
\]
for all $t\in\lbrack\alpha,\infty)_{\mathbb{T}}$.\newline Multiplying
(\ref{9}) by $e_{\ominus a}(\sigma(t),t_{0})$ and integrating the resulting
equation from $t_{0}$ to $t$ we have%
\begin{equation}
x(t)=x(t_{0})e_{a}(t,t_{0})+\int_{t_{0}}^{t}\frac{b(s)}{1+\mu(s)a(s)}%
e_{a}(t,s)x(\delta_{-}(h,s))\delta_{-}^{\Delta s}(h,s)\Delta s.\label{VAr}%
\end{equation}
Since $\delta_{-}(h,t)<\delta_{-}(\alpha,t)\leq\delta_{-}(\alpha,\alpha
)=t_{0}$ for all $t\in\lbrack t_{0},\alpha)_{\mathbb{T}}$, (\ref{zero exp})
along with Eq. (\ref{VAr}) yields%
\begin{align*}
\left\vert x(t)\right\vert  &  =e_{a}(t,t_{0})\left[  \psi(t_{0})+\int_{t_{0}%
}^{t}\frac{b(s)}{1+\mu(s)a(s)}e_{a}(t_{0},s)\psi(\delta_{-}(h,s))\delta
_{-}^{\Delta s}(h,s)\Delta s\right]  \\
&  \leq\left\Vert \psi\right\Vert \left[  e^{\int_{t_{0}}^{t}a(s)ds}%
+M\int_{t_{0}}^{t}\left\vert \frac{b(s)}{1+\mu(s)a(s)}\right\vert e^{\int
_{s}^{t}a(u)\Delta u}\Delta s\right]  \\
&  \leq\left\Vert \psi\right\Vert e^{\int_{t_{0}}^{t}a(s)ds}\left[
1+M\int_{t_{0}}^{t}\left\vert \frac{b(s)}{1+\mu(s)a(s)}\right\vert
e^{-\int_{t_{0}}^{s}a(u)\Delta u}\Delta s\right]  .
\end{align*}
The proof is complete.
\end{proof}

Notice that Theorem \ref{thm2.1} does not work for the time scales in which
\[
(t_{0},h)_{\mathbb{T}}=\varnothing.
\]
For instance, let $\mathbb{T=Z}$, $t_{0}=0$, $\delta_{-}(h,t)=t-h$ and $h=1$.
It is obvious that $(t_{0},h)_{\mathbb{Z}}=(0,1)_{\mathbb{Z}}=\varnothing$.
That is, there is no $\alpha$ so that (\ref{2.8a}) and (\ref{2.8}) hold. In
preparation for the proof of the next theorem we give the following lemma.

\begin{lemma}
Let $\mathbb{T}$ be a time scale and $t_{0}$ a fixed point. Suppose that the
shift operators $\delta_{\pm}(h,t)$ associated with the initial point $t_{0}$
are defined on $\mathbb{T}$. Suppose also that there is a delay function
$\delta_{-}(h,t)$ defined on $\mathbb{T}$. If $(t_{0},h)_{\mathbb{T}%
}=\varnothing$, then the time scale $\mathbb{T}$ is isolated (i.e.,
$\mathbb{T}$ consists only of right scattered points). Moreover,%
\begin{equation}
\sigma(t)=\delta_{+}(h,t) \label{sigma+}%
\end{equation}
for all $t\in\lbrack\delta_{-}(h,t_{0}),\infty)_{\mathbb{T}}$ or equivalently%
\begin{equation}
\sigma(\delta_{-}(h,t))=t \label{sigma-}%
\end{equation}
for all $t\in\lbrack t_{0},\infty)_{\mathbb{T}}$.
\end{lemma}

\begin{proof}
Suppose that $(t_{0},h)_{\mathbb{T}}=\varnothing$. Define $\delta_{+}%
^{0}(h,t_{0})=t_{0}$ and $\delta_{+}^{k}(h,t_{0})=\delta_{+}(h,\delta
_{+}^{k-1}(h,t_{0}))$ for $k\in\mathbb{Z}_{+}$. Since $\delta_{+}(h,t)$ is
surjective and strictly increasing we have
\[
\left(  \delta_{+}^{k-1}(h,t_{0}),\delta_{+}^{k}(h,t_{0})\right)
_{\mathbb{T}}=\delta_{-}\left(  h,\left(  \delta_{+}^{k-2}(h,t_{0}),\delta
_{+}^{k-1}(h,t_{0})\right)  _{\mathbb{T}}\right)  ,\ \ \text{for
}k=2,3,...\text{.}%
\]
Thus, one can show by induction that
\begin{equation}
\left(  \delta_{+}^{k-1}(h,t_{0}),\delta_{+}^{k}(h,t_{0})\right)
_{\mathbb{T}}=\varnothing\text{ for all }k\in\mathbb{Z}_{+}\text{.}
\label{empty}%
\end{equation}
That is, $\sigma\left(  \delta_{+}^{k-1}(h,t_{0})\right)  =\delta_{+}%
^{k}(h,t_{0})$ for $k\in\mathbb{Z}_{+}$. On the other hand, we can write%
\[
\lbrack t_{0},\infty)_{\mathbb{T}}=\cup_{k=1}^{\infty}[\delta_{+}%
^{k-1}(h,t_{0}),\delta_{+}^{k}(h,t_{0}))_{\mathbb{T}}\text{.}%
\]
Hence, for any $t\in\lbrack t_{0},\infty)_{\mathbb{T}}$ there is a $k_{0}%
\in\mathbb{Z}_{+}$ so that $t\in$ $[\delta_{+}^{k_{0}-1}(h,t_{0}),\delta
_{+}^{k_{0}}(h,t_{0}))_{\mathbb{T}}$. By (\ref{empty}) we have $t=\delta
_{+}^{k_{0}-1}(h,t_{0})$. This shows that
\[
\sigma(t)=\sigma(\delta_{+}^{k_{0}-1}(h,t_{0}))=\delta_{+}^{k_{0}}%
(h,t_{0})=\delta_{+}(h,\delta_{+}^{k_{0}-1}(h,t_{0}))=\delta_{+}(h,t)
\]
for all $t\in\lbrack\delta_{-}(h,t_{0}),\infty)_{\mathbb{T}}$. This along with
$\sigma(\delta_{-}(h,t))=\delta_{-}(h,\sigma(t))$ yields (\ref{sigma-}). The
proof is complete.
\end{proof}

\begin{theorem}
\label{thm2.1a} Let $a\in\mathcal{R}^{+}$, $Q\in\mathcal{R}$. Assume the
hypothesis of Lemma \ref{lem2}. If $(t_{0},h)_{\mathbb{T}}=\varnothing$, then
any solution $x(t)=x(t,t_{0},\varphi)$ of (\ref{9}) satisfies the exponential
inequality%
\[
\left\vert x(t)\right\vert \leq\sqrt{\left(  1+\frac{1}{\lambda}\right)
V(t_{0})}e^{\frac{1}{2}\int_{t_{0}}^{t}Q(s)\Delta s}%
\]
for all $t\in\lbrack t_{0},\infty)_{\mathbb{T}}$.
\end{theorem}

\begin{proof}
Let $H$ be defined by
\begin{equation}
H(t)=%
{\displaystyle\int\limits_{\delta_{-}(h,t)}^{t}}
\Delta s%
{\displaystyle\int\limits_{s}^{t}}
b(\delta_{+}(h,u))^{2}x(u)^{2}\Delta u. \label{H}%
\end{equation}
From (\ref{8a}), (\ref{sigma-}), and (\ref{holder}) we get%
\begin{align*}
H(t)  &  =\int_{\delta_{-}(h,t)}^{t}\left[  \sigma(u)-\delta_{-}(h,t)\right]
b(\delta_{+}(h,u))^{2}x(u)^{2}\Delta u\\
&  \geq\left[  \sigma(\delta_{-}(h,t))-\delta_{-}(h,t)\right]  \int
_{\delta_{-}(\alpha,t)}^{t}b(\delta_{+}(h,u))^{2}x(u)^{2}\Delta u\\
&  =\beta(t)\int_{\delta_{-}(\alpha,t)}^{t}b(\delta_{+}(h,u))^{2}%
x(u)^{2}\Delta u\\
&  \geq\left(
{\displaystyle\int\limits_{\delta_{-}(h,t)}^{t}}
b(\delta_{+}(h,u))x(u)\Delta u\right)  ^{2}\text{.}%
\end{align*}
Hence, by (\ref{12}) we have
\begin{align*}
V(t)  &  =A^{2}(t)+\lambda H(t)\\
&  \geq\left(  x(t)+%
{\displaystyle\int\limits_{\delta_{-}(h,t)}^{t}}
b(\delta_{+}(h,s))x(s)\Delta s\right)  ^{2}\\
&  +\lambda\left(
{\displaystyle\int\limits_{\delta_{-}(h,t)}^{t}}
b(\delta_{+}(h,u))x(u)\Delta u\right)  ^{2}\\
&  =\left(  1-\frac{1}{1+\lambda}\right)  x^{2}(t)\\
&  +\left[  \frac{1}{\sqrt{1+\lambda}}x(t)+\sqrt{1+\lambda}\left(
{\displaystyle\int\limits_{\delta_{-}(h,t)}^{t}}
b(\delta_{+}(h,u))x(u)\Delta u\right)  \right]  ^{2}\\
&  \geq\left(  1-\frac{1}{1+\lambda}\right)  x^{2}(t).
\end{align*}
This along with (\ref{v6}) yields%
\[
\left\vert x(t)\right\vert \leq\sqrt{\left(  1+\frac{1}{\lambda}\right)
V(t_{0})}e^{\frac{1}{2}\int_{t_{0}}^{t}Q(s)\Delta s}.
\]
The proof is complete.
\end{proof}

In the next corollary, we summarize the results obtained in Theorem
\ref{thm2.1} and Theorem \ref{thm2.1a}.

\begin{corollary}
\label{rem stability} Assume the hypothesis of Lemma \ref{lem2}. Let
$a\in\mathcal{R}^{+}$ and $Q\in\mathcal{R}$. Suppose that there exists a
$\lambda>0$ such that (\ref{11}) holds for all $t\in\lbrack t_{0}%
,\infty)_{\mathbb{T}}$.

\begin{enumerate}
\item If there exists an $\alpha\in(t_{0},h)_{\mathbb{T}}$ such that
(\ref{2.8a}) and (\ref{2.8}) hold, then any solution $x(t)=x(t,t_{0},\varphi)$
of (\ref{9}) satisfies
\[
\left\vert x(t)\right\vert \leq\sqrt{\frac{2}{\left(  1-\frac{1}{\xi
(t)}\right)  }V(t_{0})}e^{-\frac{1}{2}\int_{t_{0}}^{\delta_{-}(\alpha
,t)}\left[  \lambda\left(  \beta(s)+\mu(s)\right)  b(\delta_{+}(h,s))^{2}%
+\mu(s)Q^{2}(s)\right]  \Delta s}%
\]
Thus, if%
\[
\lim\limits_{t\rightarrow\infty}\int_{t_{0}}^{\delta_{-}(\alpha,t)}\left[
\lambda\left(  \beta(s)+\mu(s)\right)  b(\delta_{+}(h,s))^{2}+\mu
(s)Q^{2}(s)\right]  \Delta s=\infty,
\]
then the zero solution of Eq. (\ref{9}) is exponentially stable.

\item If $(t_{0},h)_{\mathbb{T}}=\varnothing$, then any solution
$x(t)=x(t,t_{0},\varphi)$ of (\ref{9}) satisfies%
\[
\left\vert x(t)\right\vert \leq\sqrt{\left(  1+\frac{1}{\lambda}\right)
V(t_{0})}e^{-\frac{1}{2}\int_{t_{0}}^{t}\left[  \lambda\left(  \beta
(s)+\mu(s)\right)  b(\sigma(s))^{2}+\mu(s)Q^{2}(s)\right]  \Delta s}.
\]
Thus, if%
\[
\lim\limits_{t\rightarrow\infty}\int_{t_{0}}^{t}\left[  \lambda\left(
\beta(s)+\mu(s)\right)  b(\sigma(s))^{2}+\mu(s)Q^{2}(s)\right]  \Delta
s=\infty,
\]
then the zero solution of Eq. (\ref{9}) is exponentially stable.
\end{enumerate}
\end{corollary}

Let $q>1$, $\mathbb{T}=\overline{q^{\mathbb{Z}}}=\left\{  0\right\}
\cup\left\{  q^{n}:n\in\mathbb{Z}\right\}  $, $\delta_{-}(h,t)=q^{-h}t$, and
$h\in\mathbb{Z}_{+}$. Then, Eq. (\ref{9}) turns into the $q$-difference
equation%
\begin{equation}
D_{q}x(t)=a(t)x(t)+b(t)x(q^{-h}t)q^{-h},\ \ t\in\left\{  1,q,q^{2}%
,...\right\}  , \label{qdifference}%
\end{equation}
where $D_{q}x(t)=\frac{x(qt)-x(t)}{(q-1)t}$. Next, we use Corollary
\ref{rem stability} to derive a stability criteria for the $q$-difference
equation (\ref{qdifference}).

\begin{example}
Suppose that $1+\mu(t)a(t)>0$, $1+\mu(t)Q(t)\neq0$, and%
\[
-\frac{\lambda q^{-h}}{\varpi(t)+\lambda(\varpi(t)+\mu(t))}\leq Q(t)\leq
-\lambda\left(  \varpi(t)+\mu(t)\right)  b(\delta_{+}(h,t))^{2}-\mu(t)Q^{2}(t)
\]
for all $t\in\left\{  1,q,q^{2},...\right\}  $, where $\varpi(t):=t\left(
1-q^{-h}\right)  $ and $\mu(t)=t(q-1)$.

\begin{enumerate}
\item If $(1,q^{h})_{q^{\mathbb{Z}}}\neq\varnothing$, then then condition
(\ref{2.8}) holds. By Corollary \ref{rem stability}, we conclude that any
solution $x(t)=x(t,t_{0},\varphi)$ of the $q$-difference equation
(\ref{qdifference}) satisfies the exponential inequalities%
\[
\left\vert x(t)\right\vert \leq\sqrt{\frac{2}{\left(  1-\frac{1}{\xi
(t)}\right)  }V(t_{0})}\exp\left(  \frac{1}{2}%
{\displaystyle\sum\limits_{s\in\lbrack1,q^{-\alpha}t)_{q^{Z}}}}
\mu(s)Q(s)\right)
\]
for all $t\in\lbrack q^{\alpha},\infty)_{q^{\mathbb{Z}}}$ and
\begin{align*}
\left\vert x(t)\right\vert  &  \leq\left\Vert \psi\right\Vert \exp\left(
{\displaystyle\sum\limits_{s\in\lbrack1,t)_{q^{Z}}}}
\mu(s)a(s)\right) \\
&  \times\left[  1+%
{\displaystyle\sum\limits_{s\in\lbrack1,t)_{q^{Z}}}}
G(s)\exp\left(  -%
{\displaystyle\sum\limits_{u\in\lbrack1,s)_{q^{Z}}}}
\mu(u)a(u)\right)  \right]
\end{align*}
for all $t\in\lbrack1,q^{\alpha})_{q^{\mathbb{Z}}}$, where%
\[
G(s):=q^{-h}\mu(s)\left\vert \frac{b(s)}{1+\mu(s)a(s)}\right\vert .
\]
Hence, if%
\[
\lim\limits_{t\rightarrow\infty}%
{\displaystyle\sum\limits_{s\in\lbrack1,q^{-\alpha}t)_{q^{Z}}}}
s^{2}\left[  \lambda(q-q^{-h})b(q^{h}s)^{2}+(q-1)Q^{2}(s)\Delta s\right]
=\infty,
\]
then the zero solution of Eq. (\ref{qdifference}) is exponentially stable.

\item If $(1,q^{h})_{q^{\mathbb{Z}}}=\varnothing$, then $h=1$ and
\[
\left\vert x(t)\right\vert \leq\sqrt{\left(  1+\frac{1}{\lambda}\right)
V(t_{0})}\exp\left(  \frac{1}{2}%
{\displaystyle\sum\limits_{s\in\lbrack1,t)_{q^{Z}}}}
\mu(s)Q(s)\right)
\]
Hence, if%
\[
\lim\limits_{t\rightarrow\infty}%
{\displaystyle\sum\limits_{s\in\lbrack1,t)_{q^{Z}}}}
s^{2}\left[  \lambda(q-q^{-1})b(qs)^{2}+(q-1)Q^{2}(s)\Delta s\right]
=\infty,
\]
then the zero solution of Eq. (\ref{qdifference}) is exponentially stable.
\end{enumerate}
\end{example}

In the next result, we will display a Lyapunov functional that involves
$|x|^{\Delta}$. Thus, in preparation we have the following. \newline Using the
product rule $(fg)^{\Delta}=f^{\Delta}g^{\sigma}+fg^{\Delta}$and
differentiating both sides of $x^{2}(t)=\left\vert x(t)\right\vert ^{2}$ we
obtain the derivative $\left\vert x(t)\right\vert ^{\Delta}$ as follows%
\begin{equation}
\left\vert x\right\vert ^{\Delta}=\frac{x+x^{\sigma}}{\left\vert x\right\vert
+\left\vert x^{\sigma}\right\vert }x^{\Delta}\text{ for }x\neq0. \label{sigma}%
\end{equation}
So $\left\vert x\right\vert ^{\Delta}$ depends on $\frac{x(t)}{\left\vert
x(t)\right\vert }$ and $\frac{x^{\sigma}(t)}{\left\vert x^{\sigma
}(t)\right\vert }$ (i.e., signs of $x$ and $x^{\sigma}$, respectively). Given
$x:\mathbb{T}\rightarrow\mathbb{R}$, let the sets $\mathbb{T}_{x}^{+}$ and
$\mathbb{T}_{x}^{-}$ be defined by%
\begin{align*}
\mathbb{T}_{x}^{+}  &  =\left\{  t\in\mathbb{T}:x(t)x^{\sigma}(t)\geq
0\right\}  ,\\
\mathbb{T}_{x}^{-}  &  =\left\{  t\in\mathbb{T}:x(t)x^{\sigma}(t)<0\right\}  ,
\end{align*}
respectively. The set $\mathbb{T}_{x}^{-}$ consists only of right scattered
points of $\mathbb{T}$. Since the time scale $\mathbb{T}=\mathbb{R}$ has no
any right scattered points, we have $\mathbb{T}_{x}^{-}=\varnothing$. Thus for
all differentiable functions $x:\mathbb{R}\rightarrow\mathbb{R}$, the formula
(\ref{sigma}) turns into $\left\vert x\right\vert ^{\Delta}=\frac
{x}{\left\vert x\right\vert }x^{\Delta}$. However, for an arbitrary time scale
(e.g. $\mathbb{T}=\mathbb{Z}$) the set $\mathbb{T}_{x}^{-}$ may not be empty.
For simplicity, we need to have a formula for $\left\vert x\right\vert
^{\Delta}$ which does not include $x^{\sigma}$. The next result provides a
relationship between $\left\vert x\right\vert ^{\Delta}$ and $\frac
{x}{\left\vert x\right\vert }x^{\Delta}$. Its proof can be found in
\cite{raffoul&adivar}.

\begin{lemma}
\label{lemma 3.4}\cite[Lemma 5]{raffoul&adivar} Let $x\neq0$ be $\Delta
$-differentiable. Then
\begin{equation}
\left\vert x(t)\right\vert ^{\Delta}=\left\{
\begin{array}
[c]{ll}%
\frac{x(t)}{\left\vert x(t)\right\vert }x^{\Delta}(t) & \text{if }%
t\in\mathbb{T}_{x}^{+}\\
-\frac{2}{\mu(t)}\left\vert x(t)\right\vert -\frac{x(t)}{\left\vert
x(t)\right\vert }x^{\Delta}(t) & \text{if }t\in\mathbb{T}_{x}^{-}%
\end{array}
\right.  . \label{2.2.1}%
\end{equation}

\end{lemma}

\begin{theorem}
\label{thm2.3}Define a continuous function $\eta(t)\geq0$ by%
\begin{equation}
\eta(t):=\frac{e_{a}(t,t_{0})}{1+\lambda\int_{\delta_{-}(h,t)}^{t}e_{a}%
(\delta_{+}(h,s),t_{0})\Delta s}. \label{heta}%
\end{equation}
Suppose that $a\in\mathcal{R}^{+}$ and that%
\begin{equation}
\left\vert b(t)\right\vert -\lambda\eta^{\sigma}(t)\delta_{-}^{\Delta
}(h,t)\leq0 \label{2.28a}%
\end{equation}
holds for all $t\in\lbrack t_{0},\infty)_{\mathbb{T}}$. Then any solution of
Eq. (\ref{9}) satisfies the inequality%
\begin{equation}
\left\vert x(t)\right\vert \leq V(t_{0},x_{t_{0}})e_{\gamma}(t,t_{0})\text{
for all }t\in\lbrack t_{0},\infty)_{\mathbb{T}}\text{,} \label{2.28}%
\end{equation}
where
\[
V(t_{0},x_{t_{0}}):=\left\vert x(t_{0})\right\vert +\lambda\eta(t_{0}%
)\int_{\delta_{-}(h,t_{0})}^{t_{0}}\left\vert x(s)\right\vert \Delta s,
\]
$\gamma(t):=a(t)+\lambda\widetilde{M}\eta^{\sigma}(t)$, $\widetilde{M}%
=\max\left\{  1,M\right\}  $,and $M$ is as in (\ref{delta turev}).
\end{theorem}

\begin{proof}
For convenience define%
\[
\zeta(t):=1+\lambda\int_{\delta_{-}(h,t)}^{t}e_{a}(\delta_{+}(h,s),t_{0}%
)\Delta s.
\]
Then by (\ref{4})
\begin{align}
\zeta^{\Delta}(t)  &  =\lambda e_{a}(\delta_{+}(h,t),t_{0})-e_{a}%
(t,t_{0})\delta_{-}^{\Delta}(h,t)\nonumber\\
&  =\lambda e_{a}(t,t_{0})\left[  e_{a}(\delta_{+}(h,t),t)-\delta_{-}^{\Delta
}(h,t)\right]  . \label{zeta derivative}%
\end{align}
This and a differentiation of (\ref{heta}) yield%
\begin{align}
\eta^{\Delta}(t)  &  =\frac{e_{a}(t,t_{0})}{\zeta(t)}\left(  \frac
{a\zeta(t)-\zeta^{\Delta}(t)}{\zeta^{\sigma}(t)}\right) \nonumber\\
&  =\eta(t)\left(  \frac{a\zeta(t)+a\mu(t)\zeta^{\Delta}(t)-a\mu
(t)\zeta^{\Delta}(t)-\zeta^{\Delta}(t)}{\zeta(t)+\mu(t)\zeta^{\Delta}%
(t)}\right) \nonumber\\
&  =a(t)\eta(t)-\left[  \left(  1+\mu(t)a(t)\right)  \eta(t)\frac{\zeta
(t)}{\zeta^{\sigma}(t)}\right]  \frac{\zeta^{\Delta}(t)}{\zeta(t)}\nonumber\\
&  =a(t)\eta(t)-\eta^{\sigma}(t)\frac{\zeta^{\Delta}(t)}{\zeta(t)}\nonumber\\
&  =a(t)\eta(t)+\lambda\eta^{\sigma}(t)\eta(t)\delta_{-}^{\Delta}%
(h,t)-\lambda\eta^{\sigma}(t)e_{a}(\delta_{+}(h,t),t)\nonumber\\
&  \leq\eta(t)\left[  a(t)+\lambda\widetilde{M}\eta^{\sigma}(t)\right]  ,
\label{heta derivative}%
\end{align}
where we also used $\zeta^{\sigma}(t)=\zeta(t)+\mu(t)\zeta^{\Delta}(t)$ and
\[
\left(  1+\mu(t)a(t)\right)  \eta(t)\frac{\zeta(t)}{\zeta^{\sigma}(t)}%
=\eta^{\sigma}(t).
\]
Define%
\begin{equation}
V(t,x_{t}):=\left\vert x(t)\right\vert +\lambda\eta(t)\int_{\delta_{-}%
(h,t)}^{t}\left\vert x(s)\right\vert \Delta s. \label{lyapunov}%
\end{equation}
Let $t\in\mathbb{T}_{x}^{+}\cap\lbrack t_{0},\infty)_{\mathbb{T}}$. Then by
(\ref{2.2.1}) we have $\left\vert x(t)\right\vert ^{\Delta}=\frac
{x(t)}{\left\vert x(t)\right\vert }x^{\Delta}(t)$. Differentiating
(\ref{lyapunov}) and utilizing (\ref{2.28a}) and (\ref{heta derivative}) we
arrive at%
\begin{align*}
V^{\Delta}(t,x_{t})  &  =\left\vert x(t)\right\vert ^{\Delta}+\lambda
\eta^{\Delta}(t)\int_{\delta_{-}(h,t)}^{t}\left\vert x(s)\right\vert \Delta
s\\
&  +\lambda\eta^{\sigma}(t)\left[  \left\vert x(t)\right\vert -\left\vert
x(\delta_{-}(h,t))\right\vert \delta_{-}^{\Delta}(h,t)\right] \\
&  \leq\frac{x(t)}{\left\vert x(t)\right\vert }x^{\Delta}(t)+\lambda
\eta(t)\left[  a(t)+\lambda\widetilde{M}\eta^{\sigma}(t)\right]  \int
_{\delta_{-}(h,t)}^{t}\left\vert x(s)\right\vert \Delta s\\
&  +\lambda\eta^{\sigma}(t)\left[  \left\vert x(t)\right\vert -\left\vert
x(\delta_{-}(h,t))\right\vert \delta_{-}^{\Delta}(h,t)\right] \\
&  =\left(  a(t)+\lambda\widetilde{M}\eta^{\sigma}(t)\right)  \left\vert
x(t)\right\vert +\left(  \left\vert b(t)\right\vert -\lambda\delta_{-}%
^{\Delta}(h,t)\eta^{\sigma}(t)\right)  \left\vert x(\delta_{-}%
(h,t))\right\vert \\
&  +\lambda\eta(t)\left[  a(t)+\widetilde{M}\eta^{\sigma}(t)\right]
\int_{\delta_{-}(h,t)}^{t}\left\vert x(s)\right\vert \Delta s\\
&  \leq\gamma(t)V(t,x_{t}).
\end{align*}
Similarly, if $t\in\mathbb{T}_{x}^{-}\cap\lbrack t_{0},\infty)_{\mathbb{T}}$,
then $\left\vert x(t)\right\vert ^{\Delta}=-\frac{2}{\mu(t)}\left\vert
x(t)\right\vert -\frac{x(t)}{\left\vert x(t)\right\vert }x^{\Delta}(t)$ by
(\ref{2.2.1}). Hence,%
\begin{align*}
V^{\Delta}(t,x_{t})  &  \leq\left\vert x(t)\right\vert ^{\Delta}%
+\eta(t)\left[  a(t)+\lambda\widetilde{M}\eta^{\sigma}(t)\right]  \int
_{\delta_{-}(h,t)}^{t}\left\vert x(s)\right\vert \Delta s\\
&  +\lambda\eta^{\sigma}(t)\left[  \left\vert x(t)\right\vert -\left\vert
x(\delta_{-}(h,t))\right\vert \delta_{-}^{\Delta}(h,t)\right] \\
&  \leq\left(  -\frac{2}{\mu(t)}-a(t)+\lambda\widetilde{M}\eta^{\sigma
}(t)\right)  \left\vert x(t)\right\vert \\
&  +\left(  \left\vert b(t)\right\vert -\lambda\delta_{-}^{\Delta}%
(h,t)\eta^{\sigma}(t)\right)  \left\vert x(\delta_{-}(h,t))\right\vert \\
&  +\lambda\eta(t)\left[  a(t)+\widetilde{M}\eta^{\sigma}(t)\right]
\int_{\delta_{-}(h,t)}^{t}\left\vert x(s)\right\vert \Delta s\\
&  \leq\left(  a(t)+\lambda\widetilde{M}\eta^{\sigma}(t)\right)  \left\vert
x(t)\right\vert +\lambda\eta(t)\left[  a(t)+\lambda\widetilde{M}\eta^{\sigma
}(t)\right]  \int_{\delta_{-}(h,t)}^{t}\left\vert x(s)\right\vert \Delta s\\
&  =\gamma(t)V(t,x_{t}).
\end{align*}
since $1+\mu(t)a(t)>0$ implies%
\[
-\frac{2}{\mu(t)}-a(t)<a(t).
\]
Thus,%
\begin{equation}
V^{\Delta}(t,x_{t})\leq\gamma(t)V(t,x_{t})\text{ for all }t\in\lbrack
t_{0},\infty)_{\mathbb{T}}\text{.} \label{2.29}%
\end{equation}
An integration of (\ref{2.29}) and applying the fact that $V(t,x_{t}%
)\geq\left\vert x(t)\right\vert $ we arrive at the desired result.
\end{proof}

In the next section we give a criteria for instability.

\section{A criteria for instability}

\begin{theorem}
\label{thm 3.1} Suppose there exists positive constant $D$ such that%
\begin{equation}
\beta(t)<D\leq\frac{Q(t)}{b(\delta_{+}(h,t))^{2}} \label{inscond}%
\end{equation}
for all $t\in\lbrack t_{0},\infty)_{\mathbb{T}}$, where $\beta(t)$ is as
defined in (\ref{11b}). Let the function $A$ be defined by (\ref{11a}). If
\begin{equation}
V(t)=A(t)^{2}-D\int_{\delta_{-}(h,t)}^{t}b(\delta_{+}(h,s))^{2}x(s)^{2}\Delta
s, \label{v7}%
\end{equation}
then along the solutions of Eq. (\ref{9}) we have%
\begin{equation}
V^{\Delta}(t)\geq Q(t)V(t)\text{ for all }t\in\lbrack t_{0},\infty
)_{\mathbb{T}}\text{.} \label{v8}%
\end{equation}

\end{theorem}

\begin{proof}
Let $V$ be defined by (\ref{v7}). Using (\ref{A_Q}) and (\ref{holder}) we
obtain%
\begin{align*}
V^{\Delta}(t)  &  =\left[  A(t)+A(\sigma(t))\right]  A^{\Delta}(t)-Db(\delta
_{+}(h,t))^{2}x(t)^{2}\\
&  +Db(t)^{2}x(\delta_{-}(h,t))^{2}\delta_{-}^{\Delta}(h,t)\\
&  \geq\left[  2A(t)+\mu(t)Q(t)x(t)\right]  Q(t)x(t)-Db(\delta_{+}%
(h,t))^{2}x(t)^{2}\\
&  \geq2Q(t)A(t)x(t)-Db(\delta_{+}(h,t))^{2}x(t)^{2}\\
&  =Q(t)\left[  x^{2}(t)+A^{2}(t)-\left(
{\displaystyle\int\limits_{\delta_{-}(h,t)}^{t}}
b(\delta_{+}(h,s))x(s)\Delta s\right)  ^{2}\right] \\
&  -Db(\delta_{+}(h,t))^{2}x(t)^{2}\\
&  \geq Q(t)V(t)+\left[  Q(t)-Db(\delta_{+}(h,t))^{2}\right]  x(t)^{2}.
\end{align*}
This along with (\ref{inscond}) implies (\ref{v8}).
\end{proof}

To prove the next theorem we will need to use the following lemma:

\begin{lemma}
\cite[Remarks 2]{oscillation} \ \label{remark osc} If $\varphi$ is
$rd$-continuous and nonnegative, then%
\begin{equation}
1+\int_{s}^{t}\varphi(u)\Delta u\leq e_{\varphi}(t,s)\leq\exp\left\{  \int
_{s}^{t}\varphi(u)\Delta u\right\}  \text{ for all }t\geq s.\label{osc}%
\end{equation}

\end{lemma}

\begin{theorem}
Suppose all hypotheses of Theorem \ref{thm 3.1} hold. Suppose also that
$\beta(t)$ is bounded above by $\beta_{0}$ with $0<\beta_{0}<D$. Then the zero
solution of Eq. (\ref{9}) is unstable, provided that%
\[
\lim_{t\rightarrow\infty}\int_{t_{0}}^{t}b(\delta_{+}(h,s))^{2}\Delta
s=\infty.
\]

\end{theorem}

\begin{proof}
As we did in (\ref{vin}), an integration of (\ref{v8}) from $t_{0}$ to $t$
gives%
\begin{equation}
V(t)\geq V(t_{0})e_{Q}(t,t_{0})\text{ for all }t\in\lbrack t_{0}%
,\infty)_{\mathbb{T}}\text{.} \label{3.4}%
\end{equation}
Let $V(t)$ be given by (\ref{v7}). Then%
\begin{align}
V(t)  &  =x(t)^{2}+2x(t)\int_{\delta_{-}(h,t)}^{t}b(\delta_{+}(h,s))x(s)\Delta
s+\left(  \int_{\delta_{-}(h,t)}^{t}b(\delta_{+}(h,s))x(s)\Delta s\right)
^{2}\nonumber\\
&  -D\int_{\delta_{-}(h,t)}^{t}b(\delta_{+}(h,s))^{2}x(s)^{2}\Delta s.
\label{3.5}%
\end{align}
Let $C:=D-\beta_{0}$. Then from%
\[
\left(  \sqrt{\frac{\beta_{0}}{C}}K-\sqrt{\frac{C}{\beta_{0}}}L\right)
^{2}\geq0,
\]
we have%
\[
2KL\leq\frac{\beta_{0}}{C}K^{2}+\frac{C}{\beta_{0}}L^{2}.
\]
With this in mind we arrive at%
\begin{align*}
2\left\vert x(t)\right\vert \int_{\delta_{-}(h,t)}^{t}\left\vert b(\delta
_{+}(h,s))\right\vert \left\vert x(s)\right\vert \Delta s  &  \leq\frac
{\beta_{0}}{C}x^{2}(t)\\
&  +\frac{C}{\beta_{0}}\left(  \int_{\delta_{-}(h,t)}^{t}b(\delta
_{+}(h,s))x(s)\Delta s\right)  ^{2}.
\end{align*}
A substitution of the above inequality into (\ref{3.5}) yields%
\begin{align*}
V(t)  &  \leq\left(  1+\frac{\beta_{0}}{C}\right)  x(t)^{2}+(1+\frac{C}%
{\beta_{0}})\left(  \int_{\delta_{-}(h,t)}^{t}b(\delta_{+}(h,s))x(s)\Delta
s\right)  ^{2}\\
&  -D\int_{\delta_{-}(h,t)}^{t}b(\delta_{+}(h,s))^{2}x(s)^{2}\Delta s\\
&  =\frac{D}{C}x(t)^{2}+\frac{D}{\beta_{0}}\left(  \int_{\delta_{-}(h,t)}%
^{t}b(\delta_{+}(h,s))x(s)\Delta s\right)  ^{2}\\
&  -D\int_{\delta_{-}(h,t)}^{t}b(\delta_{+}(h,s))^{2}x(s)^{2}\Delta s.
\end{align*}
Using (\ref{holder}) we find%
\[
V(t)\leq\frac{D}{C}x(t)^{2}.
\]
By (\ref{inscond}), (\ref{osc}), and (\ref{3.4}) we get%
\begin{align*}
\left\vert x(t)\right\vert  &  \geq\sqrt{\frac{C}{D}V(t_{0})e_{Q}(t,t_{0})}\\
&  \geq\sqrt{\frac{C}{D}V(t_{0})\left(  1+\int_{t_{0}}^{t}Q(s)\Delta s\right)
}\\
&  \geq\sqrt{CV(t_{0})\left(  \int_{t_{0}}^{t}b(\delta_{+}(h,s))^{2}\Delta
s\right)  }.
\end{align*}
This completes the proof.
\end{proof}

We end this paper by comparing our results to the existing ones.

\section{Some applications}

In \cite{mayr}, by means of Lyapunov's direct method the authors investigated
the stability analysis of the delay dynamic equation%
\begin{equation}
x^{\Delta}(t)=a(t)x(t)+b(t)x(\delta(t))\delta^{\Delta}(t), \label{eq1.1}%
\end{equation}
where $a:\mathbb{T\rightarrow R}$ and $b:\mathbb{T\rightarrow R}$ are
functions and $a\in\mathcal{R}^{+}$. Moreover, the delay function
$\delta:[t_{0},\infty)_{\mathbb{T}}\rightarrow\lbrack\delta(t_{0}%
),\infty)_{\mathbb{T}}$ is surjective, strictly increasing and is supposed to
have the following properties%
\[
\delta(t)<t,\ \ \ \ \delta^{\Delta}(t)<\infty,\ \ \ \delta\circ\sigma
=\sigma\circ\delta\text{.}%
\]
It is concluded in \cite[Theorem 6]{mayr} that%
\begin{equation}
\left\vert b(t)\right\vert \leq N\text{ and }a(t)<-N \label{sta}%
\end{equation}
are the sufficient conditions guaranteeing stability of the zero solution of
Eq. (\ref{eq1.1}). Next, we furnish an example to show that Theorem
\ref{thm2.1} allows us to relax condition (\ref{sta}) that leads to
exponential stability of zero solution Eq. (\ref{eq1.1}).

\begin{example}
Let $\mathbb{T=R}$, $a(t)=1$, $b(t)=-\frac{3}{2}$, $\delta(t)=t-\frac{1}{3}$,
and $N=1$. It is obvious that the condition (\ref{sta}) does not hold. So,
\cite[Theorem 6]{mayr} does not imply the stability of the zero solution of
the delayed differential equation
\begin{equation}
x^{\prime}(t)=x(t)-\frac{3}{2}x(t-\frac{1}{3}).\label{sta 1}%
\end{equation}
On the other hand, setting $\mathbb{T=R}$, $\lambda=\frac{1}{3}$, and
$\delta_{-}(h,t)=$ $t-\frac{1}{3}$ Eq. (\ref{9}) turns into (\ref{sta 1}) and
the condition (\ref{11}) becomes%
\[
-\frac{3}{4}\leq Q(t)\leq-\frac{1}{9}b(\delta_{+}(h,t))^{2},
\]
which holds for all $t\in\lbrack0,\infty)$ since $Q(t)=a(t)+b(\delta
_{+}(h,t))=-\frac{1}{2}$. One may easily verify that condition (\ref{2.8}) is
satisfied for $\delta_{-}(\alpha,t)=t-\frac{1}{6}$ and $\delta_{-}(h,t)=$
$t-\frac{1}{3}$. Thus, we conclude the exponential stability of the zero
solution of (\ref{sta 1}) by Corollary \ref{rem stability}.
\end{example}

Now, let us consider the equation%
\begin{equation}
x^{\Delta}(t)=b(t)x(\delta_{-}(h,t))\delta_{-}^{\Delta}(h,t),\ \ \ t\in\lbrack
t_{0},\infty)_{\mathbb{T}}. \label{sta3}%
\end{equation}
We observe the following by combining Corollary \ref{rem stability} and
Theorem \ref{thm2.3}.

\begin{remark}
\label{rem sta 2} Let $b\in\mathcal{R}$. Suppose that there exists a
$\lambda>0$ such that
\begin{equation}
-\frac{\lambda\delta_{-}^{\Delta}(h,t)}{\beta(t)+\lambda\left[  \beta
(t)+\mu(t)\right]  }\leq b(\delta_{+}(h,t))\leq-b(\delta_{+}(h,t))^{2}\left[
\lambda\beta(t)+(1+\lambda)\mu(t)\right]  ,\label{sta5}%
\end{equation}
holds for all $t\in\lbrack t_{0},\infty)_{\mathbb{T}}$.

\begin{enumerate}
\item If there exists an $\alpha\in(t_{0},h)_{\mathbb{T}}$ such that
(\ref{2.8a}) and (\ref{2.8}) hold and if
\begin{equation}
\lim\limits_{t\rightarrow\infty}\int_{t_{0}}^{\delta_{-}(\alpha,t)}\left[
\lambda\beta(s)+(1+\lambda)\mu(s)\right]  b(\delta_{+}(h,s))^{2}\Delta
s=\infty, \label{sta 6}%
\end{equation}
then the zero solution of Eq. (\ref{sta3}) is exponentially stable.

\item If $(t_{0},h)_{\mathbb{T}}=\varnothing$ and if
\[
\lim\limits_{t\rightarrow\infty}\int_{t_{0}}^{t}\left[  \lambda\beta
(s)+(1+\lambda)\mu(s)\right]  b(\sigma(s))^{2}\Delta s=\infty,
\]
then the zero solution of Eq. (\ref{sta3}) is exponentially stable.

\item Suppose that $a\in\mathcal{R}^{+}$ and that%
\[
\left\vert b(t)\right\vert -\lambda\eta^{\sigma}(t)\delta_{-}^{\Delta
}(h,t)\leq0
\]
holds for all $t\in\lbrack t_{0},\infty)_{\mathbb{T}}$, where%
\[
\eta(t):=\frac{1}{1+\lambda\beta(t)}.
\]
Then any solution of Eq. (\ref{9}) satisfies the inequality%
\[
\left\vert x(t)\right\vert \leq V(t_{0},x_{t_{0}})e^{\frac{1}{2}\int_{t_{0}%
}^{t}\gamma(s)\Delta s}e_{\gamma}(t,t_{0})\text{ for all }t\in\lbrack
t_{0},\infty)_{\mathbb{T}}\text{,}%
\]
where
\[
V(t_{0},x_{t_{0}}):=\left\vert x(t_{0})\right\vert +\lambda\eta(t_{0}%
)\int_{\delta_{-}(h,t_{0})}^{t_{0}}\left\vert x(s)\right\vert \Delta s,
\]
$\gamma(t):=\lambda\widetilde{M}\eta^{\sigma}(t)$, $\widetilde{M}=\max\left\{
1,M\right\}  $,and $M$ is as in (\ref{delta turev}).
\end{enumerate}
\end{remark}

In \cite[Theorem 7]{mayr}, the authors utilized fixed point theory and deduced
that the conditions
\begin{equation}
p(t):=b(\delta_{+}(h,t))\neq0\text{ for all }t\in\lbrack t_{0},\infty
)_{\mathbb{T}}, \label{stap}%
\end{equation}%
\[
\lim\limits_{t\rightarrow\infty}e_{p}(t,t_{0})=0,
\]
and%
\begin{equation}
\int_{\delta_{-}(h,t)}^{t}\left\vert p(s)\right\vert \Delta s+\int_{t_{0}}%
^{t}\left\vert \ominus p(s)\right\vert e_{p}(t,s)\left(  \int_{\delta
_{-}(h,s)}^{s}\left\vert p(u)\right\vert \Delta u\right)  \Delta s\leq N<1
\label{sta7}%
\end{equation}
lead to stability of solution $x(t,t_{0};\psi)$ of Eq. (\ref{sta3}). Notice
that \cite{mayr} generalizes all the results of \cite{raffoul}.

Moreover, Wang (see \cite[Corollary 1]{wang}) proposed the inequality%
\begin{equation}
-\frac{1}{2h}\leq a(t)+b(t+h)\leq-hb^{2}(t+h)\label{sta7-a}%
\end{equation}
as sufficient condition for uniform asymptotic stability of the zero solution
of the delay differential equation%
\[
x^{\prime}(t)=a(t)+b(t)x(t-h),\ \ h>0\text{.}%
\]

It can be easily seen that the conditions (\ref{sta7}-\ref{sta7-a}) are not
satisfied for the data given in the following example.

\begin{example}
Let $a(t)=0$, $\mathbb{T}=\mathbb{R}$, $\delta_{-}(h,t)=t-h$, and $p<0$ be
fixed. Then Eq. (\ref{9}) becomes%
\[
x^{\prime}(t)=b(t)x(t-h).
\]
We can simplify condition (\ref{sta7}) as follows%
\begin{equation}
h\left\vert p\right\vert (2-e^{pt})\leq N<1. \label{sta8}%
\end{equation}
If $h=\frac{2}{3}$, and $b(t)=-\frac{9}{10}$, then (\ref{stap}) implies
\[
h\left\vert p\right\vert (2-e^{pt})=\frac{3}{5}\left(  2-e^{-\frac{9}{10}%
t}\right)  \geq1
\]
for all $t\geq-\frac{10}{9}\ln\left(  \frac{1}{3}\right)  \cong1.22$. Thus,
the condition (\ref{sta8}) does not hold. On the other hand, for $h=\frac
{2}{3}$ and $\lambda=\frac{3}{2}$, condition (\ref{sta5}) turns into%
\[
-\frac{9}{10}\leq b(\delta_{+}(h,t))\leq-b(\delta_{+}(h,t))^{2}.
\]
The last inequality holds for $b(t)=-\frac{9}{10}$. In addition, setting
$\delta_{-}(\alpha,t)=t-\frac{1}{3}$ one may easily verify that conditions
(\ref{2.8a}), (\ref{2.8}), and (\ref{sta 6}) are satisfied. Hence, the first
part of Remark \ref{rem sta 2} yields exponential stability while
\cite[Theorem 7]{mayr} and \cite[Corollary 1]{wang} cannot.
\end{example}

\end{document}